\providecommand{\linearexample}{Pictures/heatmap_picture/}

\documentclass[%
  onecolumn
   , hidempi
]{mpi2015-cscpreprint}

\usepackage[american]{babel}

\usepackage{graphicx}

\usepackage{amssymb}
\usepackage{amsthm}
\usepackage{amsmath}

\usepackage{caption}
\usepackage{subcaption}

\usepackage{graphics} 
\usepackage{epsfig} 
\usepackage{tikz}
\usepackage{xcolor}
\usepackage[latin1]{inputenc}
\usepackage[T1]{fontenc}

\usepackage{cleveref}
\usepackage[software, hardware]{mymacros}

\usepackage{multirow}

\usepackage{comment}

\usepackage{longtable,array,amsmath,booktabs,tabularx,ragged2e}

\usetikzlibrary{arrows,automata, quotes, calc, trees, positioning, chains, shapes.geometric, decorations.pathreplacing, decorations.pathmorphing, shapes, matrix, shapes.symbols}

\usepackage{circuitikz}

\usetikzlibrary{arrows}
\usepackage{algorithm}

\usepackage[noend]{algpseudocode}
\usepackage{neuralnetwork}

\usepackage{lscape}
\usepackage[colorinlistoftodos,prependcaption,textsize=footnotesize]{todonotes}

\usepackage{hyperref}
\usepackage{enumitem}

\newcommand{\sindy}{\texttt{SINDy}}
\newcommand{\DNN}{\texttt{DNN}}
\newcommand{\RK}{\texttt{RK}}
\newcommand{\ineuralsindy}{\texttt{iNeural-SINDy}}
\newcommand{\deepymod}{\texttt{DeePyMoD}}
\newcommand{\rksindy}{\texttt{RK4-SINDy}}
\newcommand{\weaksindy}{\texttt{Weak-SINDy}}
\newcommand{\scenea}{\texttt{Scene\_A}}
\newcommand{\sceneb}{\texttt{Scene\_B}}

\begin{document}
  

\title{A Robust SINDy Approach by Combining Neural Networks and an Integral Form}

\author[$\ast$]{Ali Forootani}
\affil[$\ast$]{ Max Planck Institute for Dynamics of Complex Technical Systems, 39106 Magdeburg, Germany.\authorcr
  \email{forootani@mpi-magdeburg.mpg.de}, \orcid{0000-0001-7612-4016}}
  
\author[$\dagger$]{Pawan Goyal}
\affil[$\dagger$]{Max Planck Institute for Dynamics of Complex Technical Systems, 39106 Magdeburg, Germany.\authorcr
  \email{goyalp@mpi-magdeburg.mpg.de}, \orcid{0000-0003-3072-7780}}

\author[$\dagger\dagger$]{Peter Benner}
\affil[$\dagger\dagger$]{Max Planck Institute for Dynamics of Complex Technical Systems, 39106 Magdeburg, Germany.\authorcr
	\email{benner@mpi-magdeburg.mpg.de}, \orcid{0000-0003-3362-4103}}

\shorttitle{iNeural-SINDy--A Robust SINDy Approach}
\shortauthor{A. Forootani, P. Goyal, P. Benner}
\shortdate{}
  
\keywords{Spare regression, discovering governing equations, neural networks, nonlinear system identification, Runge-Kutta scheme.}

  
\abstract{%
The discovery of governing equations from data has been an active field of research for decades. One widely used methodology for this purpose is sparse regression for nonlinear dynamics, known as \sindy. Despite several attempts, noisy and scarce data still pose a severe challenge to the success of the \sindy~approach. In this work, we discuss a robust method to discover nonlinear governing equations from noisy and scarce data. To do this, we make use of neural networks to learn an implicit representation based on measurement data so that not only it produces the output in the vicinity of the measurements but also the time-evolution of output can be described by a dynamical system. Additionally, we learn such a dynamic system in the spirit of the \sindy~framework. Leveraging the implicit representation using neural networks, we obtain the derivative information---required for \sindy---using an automatic differentiation tool. To enhance the robustness of our methodology, we further incorporate an integral condition on the output of the implicit networks. Furthermore, we extend our methodology to handle data collected from multiple initial conditions. We demonstrate the efficiency of the proposed methodology to discover governing equations under noisy and scarce data regimes by means of several examples and compare its performance with existing methods.
}

 \novelty{\begin{itemize}
		\item We study the problem of discovering governing equations using noisy and scarce data through the lens of \sindy.
		\item We utilize neural network capabilities to avoid the requirement of explicit derivative information.
		\item In most scenarios, we observe an improved performance of the proposed methodology. 
\end{itemize}}

\maketitle

  
\section{Introduction}\label{sec:intro}

System identification is a crucial aspect of understanding and modeling the dynamics of various physical, chemical, and biological systems. Over the years, various powerful and efficient system identification techniques have been developed, and these methods have been applied in a wide range of applications, see, e.g., \cite{lennart1999system,VanM96,Tangirala_2018}. Traditionally, system identification techniques rely on prior model hypotheses. With a linear model hypothesis, several methodologies have been proposed; see, e.g.,~\cite{lennart1999system,VanM96}. However, for nonlinear system identification, defining a prior is challenging, and it is often done with the help of practitioners. Despite several earlier works \cite{NarP90,SuyVdM96,crutchfield1987equations}, nonlinear system identification is still an active and exciting research field. Towards automatic nonlinear system identification, generic algorithms and symbolic regression have shown their effectiveness and promises in discovering governing nonlinear equations using measurement~\cite{schmidt2011automated,schmidt2009distilling}. However, their computational expenses remain undesirable.

Instead of building suitable functions in the spirit of symbolic regression, there has been a focus on sparsity-promoting approaches for nonlinear system identification \cite{wang2011predicting,schaeffer2013sparse,brunton2016discovering}. They rely on the assumption that nonlinear dynamics can be defined by a few nonlinear basis functions from a dictionary with a large collection of nonlinear basis functions. Such a technique enables the discovery of interpretable, parsimonious, and generalizable models that balance precision and performance. It is nowadays widely referred to as \sindy\  \cite{brunton2016discovering}.  
\sindy\ has been employed for a handful number of challenging model discovery problems such as fluid dynamics \cite{Loiseau_2018}, plasma dynamics \cite{Dam_2017}, turbulence closures \cite{Beetham_2020}, mesoscale ocean closures \cite{Zanna_2020}, nonlinear optics \cite{Sorokina_2016}, computational chemistry \cite{Boninsegna_2018}, and numerical integration \cite{Thaler_2019}.
Moreover, the results of the \sindy\ have been extended widely to many applications, such as nonlinear model predictive control \cite{Kaiser_2018}, rational functions \cite{Kaheman_2020,Goyal_2022}, enforcing known conservation laws and symmetries \cite{Loiseau_2018}, promoting stability \cite{Kaptanoglu_2021}, generalizations for stochastic dynamics \cite{Callaham_2021}, from Bayesian perspective \cite{zhang2018robust}.

Often, \sindy~approaches require a reliable estimate of the derivative information, making them very challenging for noisy and scare data regimes. Blending numerical methods \cite{Schaeffer_2017,kang2021ident,keller2021discovery,Goyal_2022,Weak_SINDy} and weak formulations of differential equations \cite{Weak_SINDy} avoid these requirements, but their performance still deteriorates for low signal-to-noise measurements. In addition, the method in \cite{Weak_SINDy} relies on the choice of basis functions that allow to write differential equations in a weak formulation. The work in \cite{Fasel_2022_ensemble} utilizes the concepts of an ensemble to improve the predictions, but it still requires reliable estimates of derivatives to some extent. To discover governing equations from noisy data, the authors in \cite{kaheman2022automatic} proposed a scheme that aims to decompose the noisy signals into clean signals and the noise using a Runge-Kutta-based integration method. However, the scheme explicitly estimates the noise, making it harder to scale, and requires all the dependent variables to be available at the same time grid. 

Recently, applications of deep neural networks (\DNN) have received attention in sparse regression model discovery methods. For instance, in \cite{Deepymod_2021}, a deep learning-based discovery algorithm has been employed to identify underlying (partial) differential equations. However, therein, only a single trajectory is considered to recover governing equations, but in many complex processes, we might require data for different parameters and initial conditions to explore rich dynamics, thus, the reliable discovery of governing equations. Furthermore, the work \cite{Deepymod_2021} discovers governing equations based on estimating derivative information using automatic differential tools. However, we know that differential equations can also be written in the integral form, whereby the numerical approaches can employed as well, see, e.g., \cite{Goyal_2022,chen2018neural}. 

In this paper, we discuss an approach, namely, \ineuralsindy, for the data-driven discovery of nonlinear dynamical systems using noisy data from the lens of \sindy. For this, we make use of \DNN~to learn an implicit representation based on the given data set so that the network outputs denoised data, which is later utilized for the sparse regression to discover governing equations. To solve the sparse regression problem, we make use of not only automatic differential tools but also integral forms for differential equations. As a result, we observe a robust discovery of governing equations. We note that such a concept has recently been used in the context of neural ODEs in \cite{goyal2022neural} to learn black-box dynamics using noisy and scarce data. We further discuss how to incorporate the data coming from multiple initial conditions.

The rest of this paper is organized as follows. \Cref{preliminaries} briefly recalls the \sindy~approach \cite{brunton2016discovering}. In \Cref{proposed_algorithm}, we propose a novel methodology for sparse regression to learn underlying governing equations by making use of \DNN, automatic differential tools, and numerical methods. Furthermore, in \Cref{sec:multiple_initials}, we discuss its extension for multiple initial conditions and different parameters. In \Cref{Simulation}, we demonstrate the proposed framework by means of various synthetic noisy measurements and present a comparison with the current state-of-the-art approaches. Finally, \Cref{Conclusion} concludes the paper with a brief summary and future research avenues.

\section{A Brief Recap of SINDy}\label{preliminaries}
The \sindy~algorithm is a nonlinear system identification approach which is based on the hypothesis that governing equations of a nonlinear system can be given by selecting a few suitable basis functions, see, e.g., \cite{brunton2016discovering}. Precisely, it aims at identifying a few basis functions from a  dictionary, containing a large number of candidate basis functions. In this regard, sparsity-promoting approaches can be employed to discover parsimonious nonlinear dynamical systems to have a good trade-off for model complexity and accuracy \cite{Brunton_2014_compressive, Mackey_2014_compressive}. Consider the problem of discovering nonlinear systems of the form:
\begin{equation}\label{ode_1}
\dot{\bx}(t) = \mathbf{f}(\bx(t)),
\end{equation}
where $\bx(t)= \left[\bx_1(t), \bx_2(t),\ldots,\bx_n(t)\right]^\top \in \mathbb{R}^n$ denotes the state at time $t$, and $\mathbf{f}(\bx) : \R^{n} \rightarrow \R^n$ is a nonlinear function of the state $\bx(t)$.

Towards discovering the function $\mathbf{f}$ in \eqref{ode_1}, which defines the vector field or dynamics of the underlying system, we start by collecting time-series data of state ${\bx(t)}$. Let us further assume to have time derivative information of the state. If it is not readily available, we can approximate it using numerical methods, e.g., a finite difference scheme. Thus, consider that the data $\{\bx(t_0), \ldots, \bx(t_\cN)\}$ and its derivative $\{\dot\bx(t_0), \ldots, \dot\bx(t_\cN)\}$ are given.
In the next step, we assemble the data in matrices as follows:
\begin{equation}\label{data_matrix}
{\bX}=
\begin{bmatrix}
{\bx(t_1)^\top}\\{\bx(t_2)^\top}\\ \mathbf{\vdots} \\ {\bx(t_\cN)^\top}
\end{bmatrix}
= \begin{bmatrix}
{\bx_1(t_1)} & {\bx_2(t_1)} & \mathbf{\cdots} & {\bx_n(t_1)} \\
{\bx_1(t_2)} & {\bx_2(t_2)} & \mathbf{\cdots} & {\bx_n(t_2)} \\
\mathbf{\vdots} & \mathbf{\vdots} & \mathbf{\ddots} & \mathbf{\vdots} \\
{\bx_1(t_\cN)} & {\bx_2(t_\cN)} & \mathbf{\cdots} & {\bx_n(t_\cN)}
\end{bmatrix},
\end{equation}
where each row represents a snapshot of the state. Similarly, we can write the time derivative as follows:
\begin{equation}\label{data_matrix_derivative}
{\dot{\bX}} =
\begin{bmatrix}
{\dot{\bx}(t_1)^\top}\\ {\dot{\bx}(t_2)^\top}\\ \mathbf{\vdots} \\ {\dot{\bx}(t_\cN)^\top}
\end{bmatrix}
= \begin{bmatrix}
{\dot{\bx}_1(t_1)} & {\dot{\bx}_2(t_1)} & \mathbf{\cdots} & {\dot{\bx}_n(t_1)} \\
{\dot{\bx}_1(t_2)} & {\dot{\bx}_2(t_2)} & \mathbf{\cdots} & {\dot{\bx}_n(t_2)} \\
\mathbf{\vdots} & \mathbf{\vdots} & \mathbf{\ddots} & \mathbf{\vdots} \\
{\dot{\bx}_1(t_\cN)} & {\dot{\bx}_2(t_\cN)} & \mathbf{\cdots} & {\dot{\bx}_n(t_\cN)}
\end{bmatrix}.
\end{equation}
The next key building block in the \sindy\ algorithm is the construction of a dictionary $\Theta(\by)$, containing candidate basis functions (e.g., constant, polynomial or trigonometric functions). For instance, our dictionary matrix can be given as follows:
\begin{equation}\label{library_function}
\Theta(\bX) = \begin{bmatrix}
\vline & \vline & \vline & \vline & \vline & \vline & \vline & \vline & \vline \\
\mathbf{1} & {\bX} & {\bX^{\texttt P_2}} & {\bX^{\texttt P_3}} & \mathbf{\cdots} & \mathbf{\sin(X)} & \mathbf{\cos(X)} & \mathbf{\sin(2X)} & \mathbf{cos(2X)}& \cdots \\
\vline & \vline & \vline & \vline & \vline & \vline & \vline & \vline & \vline
\end{bmatrix},
\end{equation}
assume $\Theta(\mathbf{X}) \in \mathbb{R}^{m \times D}$, and in the above formulations polynomial terms are denoted by ${\bX}^{\texttt P_2}$ or ${\bX}^{\texttt P_3}$; to be more descriptive ${\bX}^{\texttt P_2}$ denotes the quadratic nonlinearities of the state ${\bX}$ as follows:
\begin{equation*}
{\bX}^{\texttt P_2}= \begin{bmatrix}
{\bx^{2}_1(t_1)} & {\bx_1(t_1)\bx_2(t_1)} & \mathbf{\cdots} & {\bx^2_2(t_1)} & {\bx_2(t_1)\bx_3(t_1)}& \mathbf{\cdots} & {\bx^2_n(t_1)} \\
{\bx^{2}_1(t_2)} & {\bx_1(t_2)\bx_2(t_2)} & \mathbf{\cdots} & {\bx^2_2(t_2)} & {\bx_2(t_2)\bx_3(t_2)} & \mathbf{\cdots} & {\bx^2_n(t_2)}\\
\vdots & \vdots & \mathbf{\ddots} & \vdots &  \vdots & \ddots & \vdots \\
{\bx^{2}_1(t_\cN)} & {\bx_1(t_\cN)\bx_2(t_\cN)} & \mathbf{\cdots} & {\bx^2_2(t_\cN)} & {\bx_2(t_\cN)\bx_3(t_\cN)} & \mathbf{\cdots} & {\bx^2_n(t_\cN)}
\end{bmatrix}.
\end{equation*}
In this setting, each column of the dictionary $\Theta(\bx)$ denotes a candidate function in defining the function $\mathbf{f}(\bx)$ in \eqref{ode_1}. 
We are interested in identifying a few candidate functions from the dictionary $\Theta$ so that a weighted sum of these selected functions can describe the function $\mathbf{f}$. For this, we can set up a sparse regression formulation to achieve this goal. Precisely, we seek to identify a sparse vector $\Xi = [\mathbf{\xi}_1,\ \mathbf{\xi}_2, \mathbf{\dots}, \ \mathbf{\xi}_n]$, where $\mathbf{\xi}_i^\top \in \R^m $ with $m$ denoting the number of columns in $\Theta$, that determines which features from the dictionary are active and their corresponding coefficients.

The \sindy\ algorithm formulates the sparse regression problem as an optimization problem as follows. Given a set of observed data $\bX$ and the corresponding time derivatives $\dot{\bX}$, the goal is to find the sparsest matrix $\Xi$ that fulfills the following:
\[\dot{\bX}=\Theta(\bX) \Xi.\]
However, finding such a matrix is an NP hard problem. Therefore, there is a need to come up a sparsity promoting regularization, and in this category, LASSO is a widely known approach \cite{friedman2001elements,tibshirani1996regression}. Despite its success, it is unable to yield matrix which is the sparsest, or the approaches, e.g., discussed in \cite{beck2013sparsity,yang2016sparse}, require a prior information about how many non-zeros elements are expected in the matrix $\Xi$, which is not known. On the other hand, the authors in \cite{brunton2016discovering} discuss a sequential thresholding approach, where simple least squares problems are solved iteratively, and at each step, coefficients below a given tolerance are pruned. Analysis of such an algorithm is discussed in \cite{zhang2019convergence}. We summarize the \sindy~approach in \Cref{STLS_SINDy_alg}. Moreover, we mention that other regularization schemes or heuristics are discussed in \cite{Rudy_2017, Schaeffer_2017, Kaptanoglu_2021, Reinbold_2020, Goyal_2022}, but in this work, we focus only on the sequential thresholding approach, similar to in \Cref{STLS_SINDy_alg} due to its simplicity.
%


\alglanguage{pseudocode}
\begin{algorithm}[tb]
	\small
	\caption{\sindy\ algorithm \cite{Fasel_2022_ensemble}}
	\label{STLS_SINDy_alg}
	\hspace*{\algorithmicindent} \textbf{Input:} Dictionary $\Theta$, time-series data $\bX$, time derivative information ${\dot{\bX}}$, threshold value $\texttt{tol}$, and  maximum iterations $ {\texttt{max-iter}}$.
	 \\
	\hspace*{\algorithmicindent} \textbf{Output:} Estimated coefficients $\Xi$ that define governing equations for nonlinear systems.
	\begin{algorithmic}[1]
		\State $\Xi = (\Theta^\top \Theta ) \big \backslash \Theta^\top \dot{\bX}$ \Comment{For initial guess, solving a least-squares problem}
		\State $k = 1$
		\While{$k < {\texttt{max-iter}}$}
		\State \texttt{small\_inds} = (\textbf{abs}($\Xi$) < $\texttt{tol}$) \Comment{identifying small coefficients}
		\State $\Xi$(\texttt{small\_inds}) = 0 \Comment{excluding small coefficients}
		\State Solve $\Xi = (\Theta^\top \Theta ) \big \backslash \Theta^\top \dot{\bX}$ subject to $\Xi$(\texttt{small\_inds}) = 0 
		\State $k = k +1$
		\EndWhile
	\end{algorithmic}
	\vspace{0.2cm}
\end{algorithm}

\section{\ineuralsindy: Neural Networks and Integrating Schemes  Assisted \sindy\ Approach} \label{proposed_algorithm}

A challenge in the classical \sindy\ approach discussed in the previous section is the availability of an accurate estimate of the derivative information. If the derivative information is inaccurate, the resulting sparse model may not accurately capture the underlying system dynamics. 

In this section, we present an approach that combines \sindy\ framework with a numerical integration scheme and neural networks in a particular way so that a robust discovery of governing equations can be made amid poor signal-to-noise ratio and irregularities in data. The methodology is inspired by the work \cite{goyal2022neural}. The main components of the methodology are as follows. For given noisy data, we aim to learn an implicit representation using a neural network based on the noisy data so that the network yields denoised data but still in the vicinity of the noisy collected data, and governing equations describing the dynamics of the denoised data can be obtained by employing \sindy. For \sindy, we utilize automatic differential tools to obtain the derivative information via the network and also make use of an integral form of the differential equations. In the following, we make these discussions more precise.

Consider noisy data ${\by(t)} \in \mathbb{R}^n$ at the time instances $\{t_0,\dots,t_{\cN}\}$, i.e., $\{\by(t_0), \ldots, \by(t_\cN)\}$. Moreover, $\by(t) = \bx(t) + \epsilon(t)$, where $\bx(t)$ and $\epsilon(t)$ denote clean data and noise, respectively. Under this setting, we aim to discover the structure of vector field $\mathbf{f}$ by identifying the most active terms in the dictionary $\Theta$ so that it satisfies as follows:
\begin{equation}\label{dynamical_sys}
\dot{\bx}(t)=\mathbf{f}(\bx).
\end{equation} 
Note that we do not know $\epsilon$'s. In order to learn $\mathbf{f}$ from $\by$, we blend three ingredients together, which are discussed in the following.

\begin{enumerate}[label=(\alph*)]
	\item \textbf{Sparse regression assumption:} In our setting, we utilize the principle of \sindy, which we discussed in the previous section. This means that the system dynamics (or vector field defining dynamics) can be represented by a few suitable terms from a dictionary of candidate functions. This allows us to obtain a parsimonious representation of dynamical systems and reduces the model complexity, leading to better generalization and interpretability of the models.
	
	\item \textbf{Automatic differential to estimate derivative information:}
	As mentioned earlier, \sindy\ algorithm requires accurate derivative information for the system, which can be challenging to obtain from experiments or to estimate using numerical methods. To cope with this issue, we make use of neural networks with its automatic differentiation (AD) feature, which is a technique used to estimate derivative information. The use of a \DNN\ in combination with the \sindy\ algorithm was earlier discussed in \cite{Deepymod_2021}, where it has been shown that the discovery of nonlinear system dynamics without explicit need of accurate derivative information \cite{Deepymod_2021}.
	
	We make use of a \DNN\ to parameterize a nonlinear mapping from time $t$ to the dependent variable $\by(t)$. To that end, let us denote a \DNN~by $\cG_\theta$, where $\theta$ contains \DNN~parameters. The input to $\cG_\theta$ is time $t$, and its output is $\by(t)$, i.e., 
	$\by(t) = \cG_\theta(t)$. However, in the case of noisy measurement $\by(t)$ at time $\{t_0,\dots,t_{\cN}\}$, we expect $\cG_\theta$ to predict outputs in the proximity of $\by$, i.e., 
	\[\by(t) \approx \bx(t) =  \cG_\theta(t),\quad t = \{t_0,\dots,t_{\cN}\}.\]
	With the sparse regression hypothesis, we aim to learn a dynamical model for $\bx$, as it can be seen as a denoised version of $\by$. For this, we construct a dictionary of possible candidate functions using $\bx$, which we denote by $\Theta\big(\bX(t)\big)$. Next, we require the derivative information of $\bx$ with respect to time $t$. Since we have an implicit representation of $\bx(t)$ using a \DNN, we can employ AD to obtain the required information. Having dictionary and derivative information, we set up a sparse regression problem as follows:
	\begin{equation}\label{Automatic_diff}
	\dot{\bX}(t) = \Theta\big(\bX(t)\big)\Xi,
	\end{equation}
	where $\Xi$ is the sparsest possible matrix, which selects the most active terms from the dictionary to define dynamics. Finding the sparsest solution is computationally infeasible; we, thus, utilize the sequential thresholding approach as discussed in \Cref{STLS_SINDy_alg} with minor modifications. Instead of solving least-squares problems in Steps $1$ and $5$ in \Cref{STLS_SINDy_alg}, we have a loss function as follows:
	\begin{equation}
	\cL := \min_{\theta,\Xi}\sum_{i=0}^\cN \lambda_1\|\by(t_i) - \bx(t_i)\| + \lambda_2 \|\dot{\bx}(t) - \Theta\big(\bx(t)\big)\Xi\|,
	\end{equation}
	where $\bx(t_i) := \cG_\theta(t_i)$, and $\lambda_{1}$ and $\lambda_{2}$ are hyperparameters.
	
	\item \textbf{Numerical integration scheme:} A dynamical system is a particle or an ensemble of particles whose state varies over time and thus obeys differential equations involving time derivatives \cite{Sprott_2003_chaos}. To predict the evolution of the dynamical system, it is necessary to have an analytical solution of such equations or their integration over time through computer simulations.
	Therefore, we aim to incorporate the information contained in the form of integration of dynamical systems while discovering governing equations via sparse regression, which is expected to make the process of discovering equations robust to the noise and scarcity of data.

	When differential equations are written in an integral form, then we do not require derivative information as well; however, the resulting optimization problem involves an integral form. In this regard, one can employ the principle of Neural-ODEs \cite{NODE} to solve efficiently such optimization problems. One can also approximate the integral form using suitable integrating schemes \cite{butcher2016numerical}, and recently, fourth-order Runge-Kutta (\RK4) scheme \cite{Goyal_2022} and linear multi-step methods \cite{buckwar2006multistep} are combined with \sindy. In this work, we make use of the \RK4\ scheme to approximate an integral.	
	
	Following \cite{Goyal_2022}, our goal is to predict the state of a dynamic system ${\bx(t_{k+1})}$ at time $t = t_{k+1}$ from the state ${\bx(t_k)}$ at time $t = t_k$, where $k \in \{0,1,\dots,\cN-1 \}$. By employing the \RK4\ scheme, ${\bx(t_{k+1})}$ can be computed as a weighted sum of four components that are the product of the time-step and gradient field information $\mathbf{f}(\cdot)$ at the specific locations. These components are computed as follows:
	\begin{equation}\label{RK}
	{\bx(t_{k+1}) \approx \bx(t_k) + \frac{1}{6} h_k (\ba_1 + 2\cdot\ba_2 + 2\cdot\ba_3 + \ba_4),\ \ \ h_k=t_{k+1} - t_{k}},
	\end{equation}
	where,
	\begin{equation*}
	\ba_1= \mathbf{f}(\bx(t_k)), \ \ \ \ba_2 = \mathbf{f} \Big( \bx\big( t_k) + h_k \frac{\ba_1}{2}   \Big),\ \ \ a_3= \mathbf{f} \Big(\bx( t_k) + h_k \frac{\ba_2}{2} \Big),\ \ \ \ba_4= \mathbf{f} \Big(\bx( t_k) + h_k \ba_3 \Big).
	\end{equation*}
	
	For the sake of simplicity with a slight abuse of a notation, the right-hand side of \eqref{RK} is denoted by $\mathbf{\mathcal{F}}_{\text{\RK4}} \big(f,\bx(t_k),h_k\big)$, i.e.,
	\begin{equation}
	{\bx(t_{k+1}) = \bx(t_k + h_k) \approx \mathcal{F}_{\text{\RK 4}} \big(\mathbf{f},\bx(t_k),h_k \big)}.
	\end{equation}
	Like the \sindy\ algorithm, we collect samples from the dynamical system at time $t = \{t_0,\dots,t_\cN\}$ and define the time step as $h_k:= t_{k+1}-t_k$. 
	
	With sparse regression assumption, we can write $\mathbf{f}(\bx) = \Theta(\bx)\Xi$, where $\Theta(\bx)$ is a dictionary and $\Xi$ is a sparse matrix. Then, we can set up a sparse regression as follows. We seek to identify the sparsest matrix $\Xi$ so that the following is minimized:
	\[\sum_k\left\|\bx(t_{k+1}) - \mathcal{F}_{\RK 4}\big(\Theta(\bx)\Xi,\bx(t_k),h_k\big)\right\|.\]
	
	When the \RK4 scheme is merged with the previously discussed \DNN\ framework, we apply a one-time ahead prediction based on
	\rksindy\ to the output of our \DNN, i.e.,
	\[\bx_{\RK4}(t_{k+1}) \approx \mathcal{F}_{\RK4}\Big( \mathbf{f},\bx(t_k),h_k \Big).\]
	
\end{enumerate}

Having all these ingredients, we combine them to define a loss function to train our \DNN\ structure, as well as to discover governing equations describing underlying dynamics. To that end, we have the following loss function:
\begin{equation}\label{objective_function}
\mathbf{\mathcal{L}} = \mu_1 \mathbf{\mathcal{L}}_{\texttt{MSE}} + \mu_2 \mathbf{\mathcal{L}}_{\texttt{deri}} + \mu_3 \mathbf{\mathcal{L}}_{\texttt{RK4}}, \ \ \ \mu_1,\mu_2,\mu_3 \in [0,1],
\end{equation}
where $\mathbf{\mathcal{L}}_{\texttt{MSE}}$ is the mean square error (\texttt{MSE}) of the output of the \DNN\ $\cG_\theta$ (denoted by $\mathbf{\hat{x}}$) with respect to the collected data ${\by}$, and $\{\mu_{1}, \mu_{2}, \mu_{3}\}$ are positive constants, determining the weight of different losses in the total loss function. It is given as
\begin{equation}\label{mse_loss}
\mathbf{\mathcal{L}}_{\texttt{MSE}} = \frac{1}{\cN} \sum_{k=1}^{\cN}\Big\|{\by(t_k) - \bx(t_k)}\Big\|_2^2.
\end{equation}
It forces the \DNN\ to produce output in the vicinity of the measurements, and ${\mu_1}$ is its weight. $\mathbf{\mathcal{L}}_{\texttt{deri}}$ is inspired by the sparse regression and aims to compute the sparse coefficient matrix $\Xi$. It is computed as follows:
\begin{equation}\label{reg_loss}
\mathbf{\mathcal{L}}_{\texttt{deri}} = \frac{1}{\cN} \sum_{k=1}^{\cN} \Big\| \dot{\bx}(t_k)  -\Theta\big(\bx(t_k)\big) \Xi \Big\|_2^2,
\end{equation} 
The term $\mathbf{\mathcal{L}}_{\texttt{RK4}}$ encodes the capabilities of the vector field to predict the state at the next time step. This is the \texttt{MSE} of the output of the \RK4 scheme and the output of \DNN,\ given as follows:
\begin{equation}\label{Rk4_loss}
\mathbf{\mathcal{L}}_{\RK4} = \frac{1}{\cN-1} \sum_{k=1}^{\cN} \Big\| \frac{1}{h_k} \left(\bx(t_{k+1}) -\mathcal{F}_{\RK4}\left( \Theta\big(\bx(t_k)\big) {\Xi},\bx(t_k),h_k \right) \right)\Big\|_2^2.
\end{equation} 
It is worth highlighting that the coefficient matrix $\Xi$ will be updated alongside the weights and biases of the \DNN, and the dictionary terms are calculated by \eqref{library_function}. Furthermore, after a certain number of epoch training, we employ sequential thresholding on $\Xi$ to remove small coefficients as sketched in \Cref{STLS_SINDy_alg}, and update the remaining parameters thereafter. We summarize the procedure in \Cref{iNeuralSINDy_alg}. Additional steps in \Cref{iNeuralSINDy_alg} are as follows. We train our network for initial iterations (denoted by {\texttt{init-iter}}) without employing sequential thresholding; this helps the \DNN~to learn the underlying dynamics of the dataset. Afterward, we employ sequential thresholding every $q$ iterations. In the rest of the paper, the proposed methodology is referred to as \ineuralsindy.


\alglanguage{pseudocode}
\begin{algorithm}[tb]
	\small
	\caption{\texttt{iNeuralSINDy:} \texttt{SINDy} combined with neural network and integral scheme for nonlinear system identification.}
	\label{iNeuralSINDy_alg}
	
	\hspace*{\algorithmicindent} \textbf{Input:} Data set $\{\by(t_0),\by(t_1),\dots,\by(t_\cN)\}$,  $\texttt{tol}$ for sequential thresholding, a dictionary containing candidate functions $\Theta$, a neural network $\cG_\theta$ (parameterized by $\theta$),  initial iteration ($ {\texttt{init-iter}}$), maximum iterations ${\texttt{max-iter}}$, and parameters $\{\mu_1,\mu_2,\mu_3\}$.
	\\
	\hspace*{\algorithmicindent} \textbf{Output:} Estimated coefficients ${\Xi}$, defining governing equations. 
	\begin{algorithmic}[1]
		
		\State Initialize the \texttt{DNN} module parameters, and the coefficients $\Xi$
		\State $k=1$
		\While {$k < {\texttt{max-iter}}$}
		\State Feed time $t_i$ as an input to the \DNN~($\cG_\theta$) and predict output $\bx$.
		\State Compute the derivative information $ \dot{\bx}$ using automatic differentiation. 
		\State Compute the cost function~\eqref{objective_function}.
		\State Update the parameters of \DNN\ ($\theta$) and the coefficient ${\Xi}$ using gradient descent. 
		\For {$k\%q ==0$ \& $k >  {\texttt{init-iter}}$}  \Comment Employing sequential thresholding after $q$ iterations
		\State \texttt{small\_inds} = (\textbf{abs}($\Xi$) < $\texttt{tol}$) \Comment{identifying small coefficients}
		\State $\Xi$(\texttt{small\_inds}) = 0 \Comment{excluding small coefficients}
		\State Update the parameters of \DNN\ ($\theta$) and the coefficient ${\Xi}$ using gradient descent, 
		\Statex \qquad\qquad while  ensuring $\Xi$(\texttt{small\_inds}) remains zero. 
		\EndFor
		\State $k=k+1$.
		\EndWhile
		
	\end{algorithmic}
	\vspace{0.2cm}
\end{algorithm}


\section{Extension to Multi-trajectories Data}\label{sec:multiple_initials}
Thus far, we have presented the discovery of governing equations using a single trajectory time series data set using a single initial condition. However, for complex dynamical processes, a single trajectory is not sufficient to describe underlying dynamics completely. Therefore, it is necessary to collect data using multiple trajectories; hence, we need to adopt our proposed methodology to account for multiple trajectories. 

\begin{figure}[tb]
	\centering
	\includegraphics[width=\textwidth]{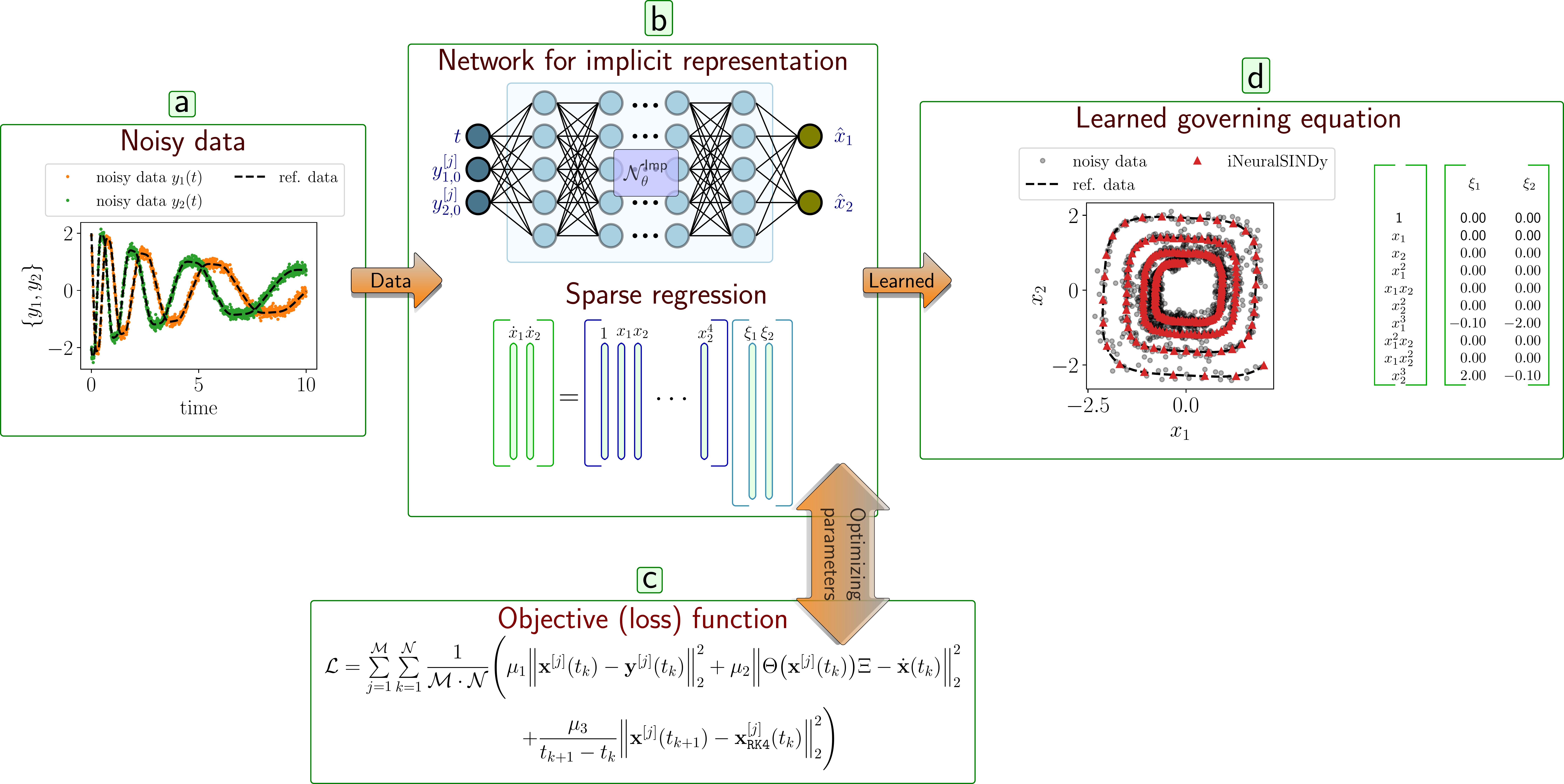}
	\caption{ A schematic diagram of the approach \ineuralsindy. (a) noisy measurement data, (b) feeding the initial condition $\big({\by_{1,0},\ \by_{2,0}}\big)$ and the time $t$ to the \texttt{DNN}, (c) using the output of the \texttt{DNN},  construct a polynomial dictionary, (d) estimating the parameters of the \texttt{DNN} and sparse vector $\Xi$  by considering a loss function.
	}
	\label{iNeuralSINDy_schematic}	
\end{figure}

To achieve this goal, we augment the input time $t$ with an initial condition so that a \DNN\ can capture the nonlinear behavior of the system with respect to different initial conditions. 
To that end, let us consider $\cM$ different trajectories with initial conditions $y_0^{[j]}$, where $j \in \{1,\ldots, \cM\}$. To reflect the multi-trajectories in our framework, we modify the architecture of the \DNN, which now takes $t_k$ and $y_0^{[j]}$ as inputs, and intend to predict $y_k^{[j]}$---that is, the state at time $t_k$ with respect to the initial condition  $y_0^{[j]}$. Then, we also adapt our loss function \eqref{objective_function} as follows:
\begin{equation}\label{objective_function_j}
\mathbf{\mathcal{L}} = \mu_1 \sum_{j=1}^{\cM} \mathbf{\mathcal{L}}_{\texttt{MSE}}^{[j]} + \mu_2 \sum_{j=1}^{\cM} \mathbf{\mathcal{L}}_{\texttt{deri}}^{[j]} + \mu_3 \sum_{j=1}^{\cM} \mathbf{\mathcal{L}}_{\texttt{RK4}}^{[j]}, \ \ \ \mu_1,\mu_2,\mu_3 \in [0,1],
\end{equation}
where
\begin{align*}
\mathbf{\mathcal{L}}_{\texttt{MSE}}^{[j]} &= \frac{1}{\cM \cdot \cN}\sum_{j=1}^{\cM} \sum_{k=1}^{\cN}\Big\|{\by^{[j]}(t_k) - \bx^{[j]}(t_k)}\Big\|_2^2,\\
\mathbf{\mathcal{L}}_{\texttt{deri}}^{[j]} &= \frac{1}{\cM \cdot \cN}
\sum_{j=1}^{\cM} \sum_{k=1}^{\cN} \Big\| \Theta\big(\bx^{[j]}(t_k)\big) \hat{\Xi} -
\dot{\bx}(t_k) \Big\|_2^2,\\
\mathbf{\mathcal{L}}_{\RK4}^{[j]} &= \frac{1}{h} \frac{1}{\cM \cdot \cN}\sum_{j=1}^{\cM} \sum_{k=1}^{\cN} \Big\| \bx^{[j]}(t_{k+1}) - \bx^{[j]}_{\RK4}(t_k) \Big\|_2^2~~\text{with}~~h = t_{k+1} - t_k.
\end{align*} 
We depict a schematic diagram of our proposed approach in \Cref{iNeuralSINDy_schematic} for such a case.

\section{Numerical Experiments} \label{Simulation}
In this section, we demonstrate the proposed methodology, the so-called \ineuralsindy, by means of several numerical examples and present a comparison with existing methodologies.
For the comparison, we primarily consider two approaches, namely \deepymod~\cite{Deepymod_2021}, and \rksindy~\cite{Goyal_2022}. \deepymod~utilizes only automatic-differential tools to estimate derivative information by constructing an implicit representation of the noisy data, while \rksindy~embeds a numerical integration scheme to avoid computation of derivative information. The proposed methodology \ineuralsindy~can be viewed as a combination of \deepymod~and \rksindy. For the chaotic Lorenz example, we also present a comparison with \weaksindy~\cite{Weak_SINDy}.

To quantify the performance of the considered methodologies, we define the following coefficient error measure for each state variable $\bx_i$:
\begin{equation}\label{error_criteria}
\cE({\bx_i}) = \left\|\Xi_{\bx_i}^{\texttt{truth}} - \Xi_{\bx_i}^{\texttt{est}}  \right \|_{1},
\end{equation}
where $\Xi_{\bx_i}^{\texttt{truth}}$ and ${\Xi}^{\texttt{est}}_{\bx_i}$ are, respectively, the true and estimated coefficients, corresponding to the state variable ${\bx_i}$, and $\|\cdot \|_1$ denotes the $l_1$-norm. A motivation to quantity each state variable separately is that their dynamics can be of different scales; thus, their coefficients might also be in a different order. Therefore, to better understand the quality of the discovered models, we analyze them separately. 
Furthermore, to observe the performance of the methodologies under the noisy data, which is often the case in real-world scenarios, we artificially generate noisy data by corrupting the clean data. 
For this, we use a white Gaussian noise $\cN(\mu, \sigma^2)$ with a zero mean $\mu=0$ and variance $\sigma^2$, where $\sigma$ denotes the standard deviation. The noise level in the data is controlled by $\sigma$, i.e., larger $\sigma$ implies more noise present in the data.
Additionally, since \ineuralsindy~and \deepymod~both involve neural networks, we also compare their performance sensitivity in two scenarios as follows:
\begin{itemize}
	\item \scenea: In the first scenario, we consider having a single initial condition and a fixed number of neurons in the hidden layers but vary the amount of training data and noise levels.
	\item \sceneb: In the second one, we consider having a single initial condition and a fixed number of training data but vary the number of neurons in the hidden layers and noise levels.
\end{itemize}

In addition, in the following, we further clarify common implementation and reproducibility details that are considered for all the examples. 

\paragraph{Data generation.}
We have generated the data synthetically by using \texttt{solve\_ivp} function from \texttt{scipy.integrate} package to solve a given set of differential equations and produce the data set. When an identification approach terminates based on the considered methodologies (e.g., \ineuralsindy, \deepymod, \rksindy, or \weaksindy), we multiply the dictionary $\Theta$ by estimated coefficient matrix $\Xi^{\texttt{est}}$ to obtain the discovered governing equations. We then make use of the \texttt{solve\_ivp} function from \texttt{scipy.integrate} to obtain time-evolution dynamics. 

Moreover, we perform a data-processing step before feeding to a neural network by mapping the minimum and maximum values to $-1$ and $1$, respectively.
The  hyper-parameters $\mu$'s in \eqref{objective_function_j} are set to $\mu_1=1$, $\mu_2=0.1$ and $\mu_3=0.1$ for \ineuralsindy. Note that we can drive \rksindy\ and \deepymod\ approaches by setting $\mu_3=0$ and $\mu_2=0$, respectively, in \eqref{objective_function_j}. 

\paragraph{Architecture.}
We use multi-layer perception networks with periodic activation functions, namely, \texttt{SIREN} \cite{SIREN}, to learn an implicit representation based on measurement data. The numbers of hidden layers and neurons will be discussed for each example separately.

\paragraph{Hardware.} For training neural networks and parameter estimations for discovering governing equations, we have used \nvidia RTX A4000 GPU with 16 GB RAM, and for CPU computations (e.g., for generating data), we have used a 12th Gen \intel~\coreifive-12600K processor with 32 GB RAM. 

\paragraph{Training set-up.}
We use the Adam optimizer \cite{kingma2014adam} to update the coefficient matrix $\Xi$ that is trained alongside the \DNN\ parameters. The threshold value (\texttt{tol}), learning rate of the optimizer, maximum iterations(${\texttt{max-iter}}$),  initial iterations (${\texttt{init-iter}}$), the iteration $q$ for employing sequential thresholding for \Cref{iNeuralSINDy_alg} will be mentioned for each example separately. 

However, we note that after each thresholding step in \Cref{iNeuralSINDy_alg}, we reset the learning rate $5\times 10^{-6}$ for \DNN\ parameters and $1 \times 10^{-2}$ for the coefficient matrix $\Xi^{\texttt{est}}$ except for the Lorenz example, which is explicitly mentioned in the Lorenz example.

\subsection{Two-dimensional damped oscillators}
In our first example, we consider the discovery of a two-dimensional linear oscillatory damped system using data. The dynamics of the oscillator can be given by
\begin{equation}\label{2_d_osci}
\begin{aligned}
\dot{\bx}_1(t) &= -0.1 {\bx_1(t)} + 2.0 {\bx_2(t)},\\
\dot{\bx}_2(t) &= -2.0 {\bx_1(t)} - 0.1 {\bx_2(t)}.
\end{aligned}
\end{equation}

\paragraph{Simulation setup:}
To generate the training data set, we consider three initial conditions in the range $[-2,2]$ for $\bx_1$ and $\bx_2$, and for each initial condition, we take $400$ equidistant points in the time interval $t \in [0,\ 10]$. Our \DNN\ architecture has three hidden layers, each having $32$ neurons. We set the number of epochs ${\texttt{max-iter}} = 15,000$ and threshold value $\texttt{tol}= 0.05$. The initial iteration {\texttt{init-iter}} is set to $5,000$ with the learning rate of $10^{-4}$ for the \DNN\ parameters and $10^{-3}$ for the coefficient matrix $\Xi^{\texttt{est}}$, and after $q=2,000$ iterations, we employ the sequential thresholding. Moreover, we construct a dictionary containing polynomials of degrees up to two.

\begin{figure}[!tb]
	\centering
	\includegraphics[width=0.7\textwidth, trim = 0cm 7.8cm 0cm 0cm, clip]{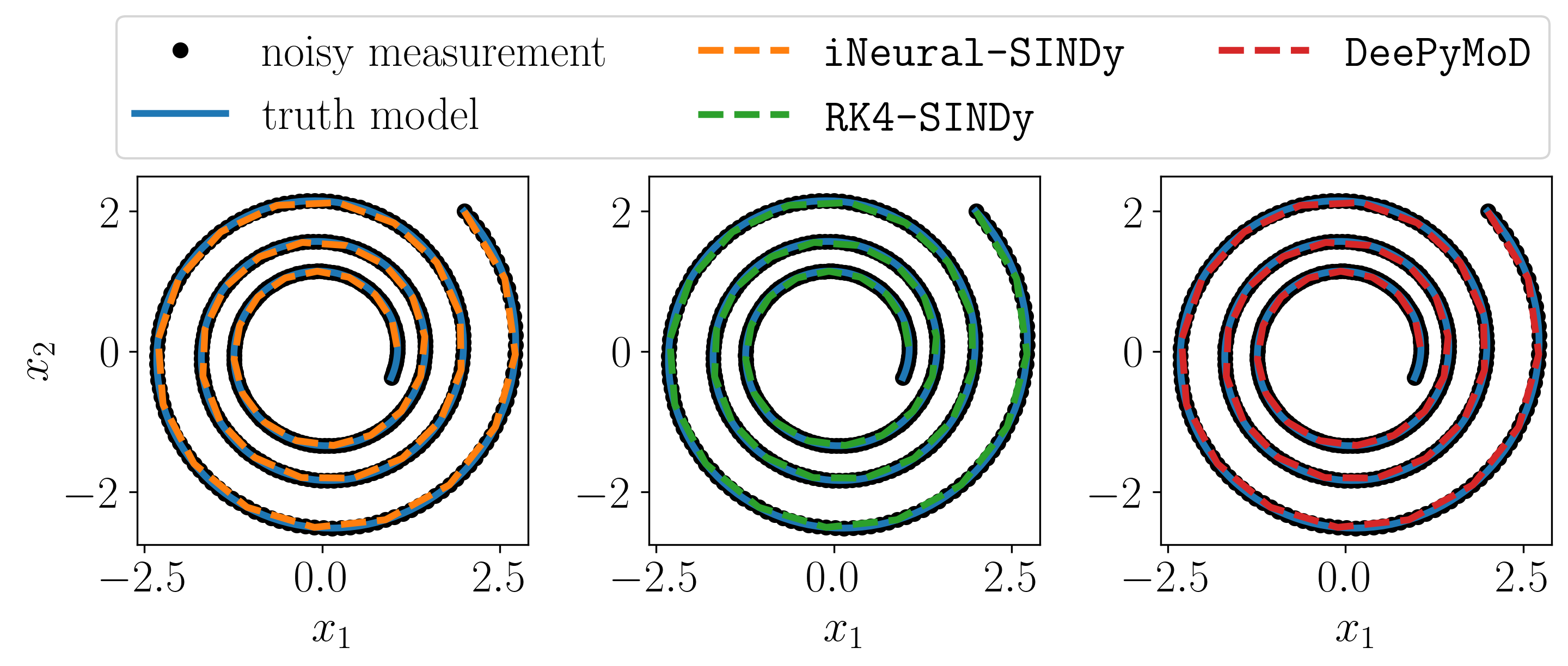} \\[2pt]
	
	\begin{subfigure}{0.49\textwidth}
		\includegraphics[width=\textwidth, trim = 0cm 0cm 0cm 2.75cm, clip]{Pictures/2_D_Oscillatory/2_D_Oscillatory_noise_0.0.png}
		\caption{noise level $\sigma = 0.00$}
		\label{2_D_Oscillatory_noise_0}
	\end{subfigure}
	\hfill
	\begin{subfigure}{0.49\textwidth}
		\includegraphics[width=\textwidth, trim = 0cm 0cm 0cm 2.75cm, clip]{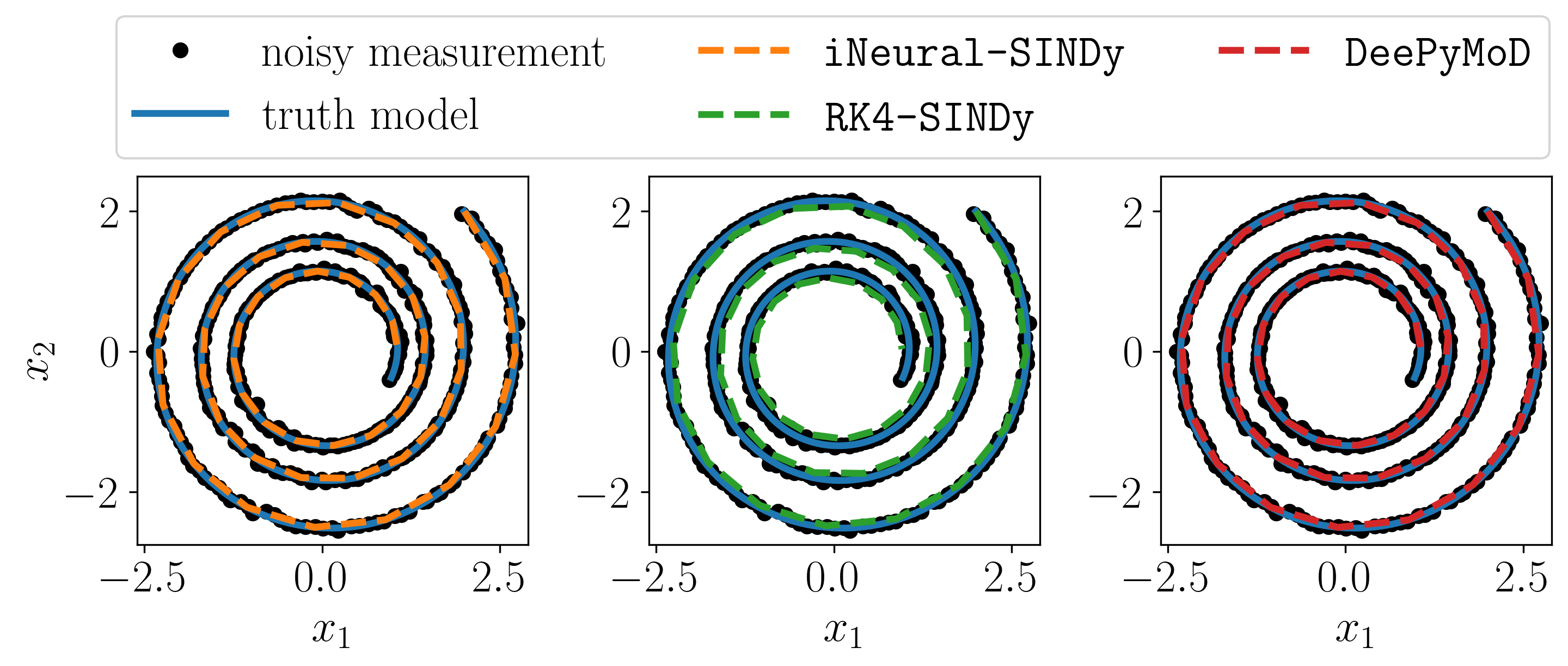}
		\caption{noise level $\sigma = 0.02$}
		\label{2_D_Oscillatory_noise_02}
	\end{subfigure}
	
	\begin{subfigure}{0.49\textwidth}
		\includegraphics[width=\textwidth, trim = 0cm 0cm 0cm 2.75cm, clip]{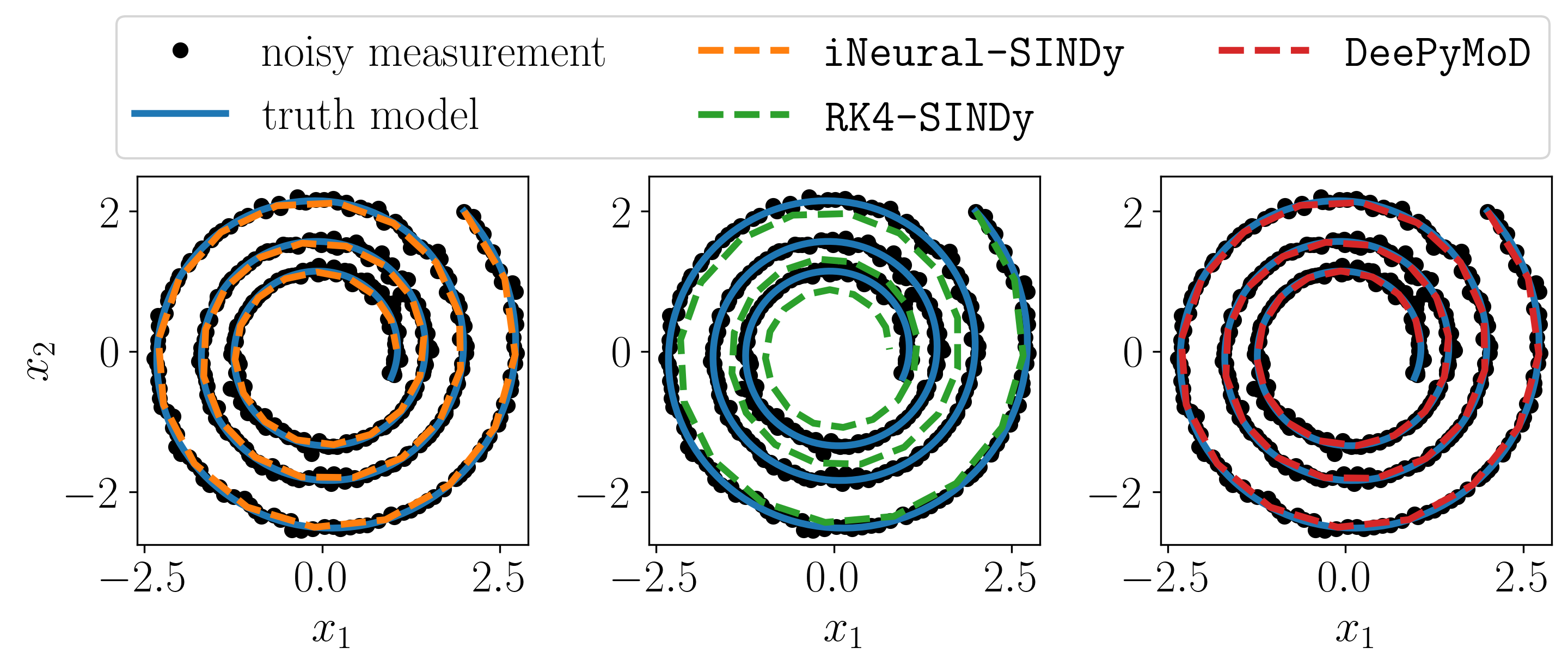}
		\caption{noise level $\sigma = 0.04$}
		\label{2_D_Oscillatory_noise_04}
	\end{subfigure}
	\hfill
	\begin{subfigure}{0.49\textwidth}
		\includegraphics[width=\textwidth, trim = 0cm 0cm 0cm 2.75cm, clip]{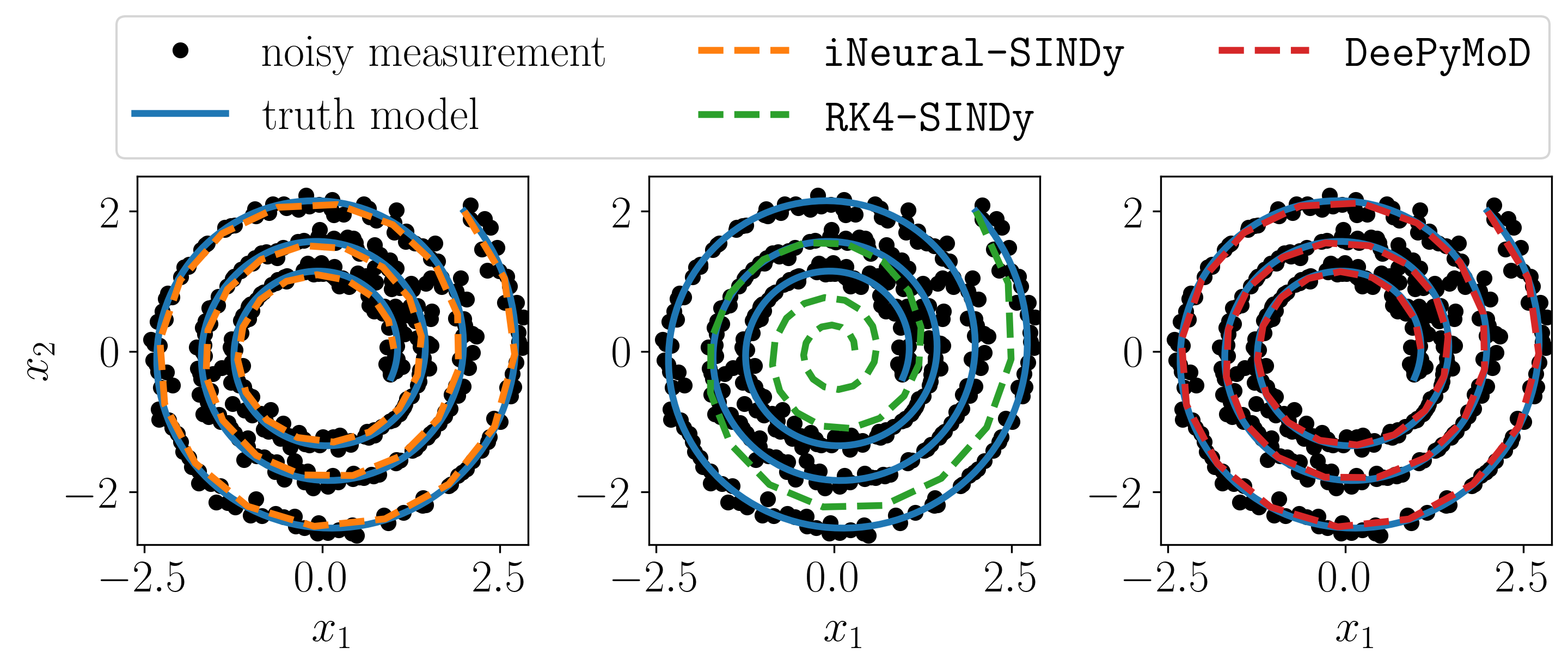}
		\caption{noise level $\sigma = 0.08$}
		\label{2_D_Oscillatory_noise_08}
	\end{subfigure}
	\caption{Linear oscillator: A comparison of the learned equations using different methods under various noise levels present in measurement with the ground truth.}
	\label{2_D_Oscilatory_noise_level}
\end{figure}

\begin{figure}[!tb]
	\centering
	\includegraphics[width=0.8\textwidth,trim = 0cm 7cm 0cm 0cm, clip]{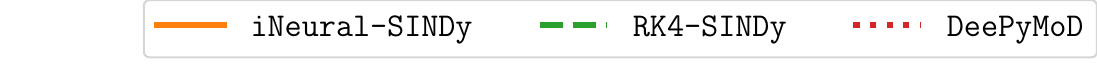}
	\begin{subfigure}{0.49\textwidth}
		\includegraphics[width=1\textwidth,trim = 0cm 0cm 0cm 1.1cm, clip]{\linearexample/2_D_Oscilator/2_D_Oscilator_noise_0.0_Coefficients__threshold_0.05coeff_1.pdf}
		\caption{Coefficients for $x_1$.}
		\label{}
	\end{subfigure}
	\begin{subfigure}{0.49\textwidth}
		\includegraphics[width=1\textwidth,trim = 0cm 0cm 0cm 1.1cm, clip]{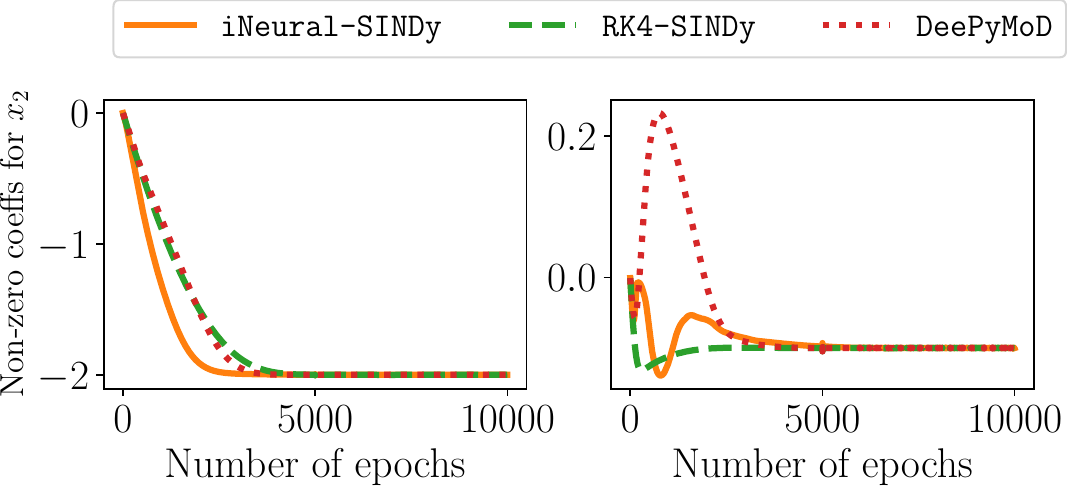}
		\caption{Coefficients for $x_2$.}
		\label{}
	\end{subfigure}
	\caption{Linear oscillator: Estimated coefficients during the training loop for \ineuralsindy, \deepymod\ and \rksindy.}
	\label{Coeff_iteration_2_D_Oscillator}	
\end{figure}

\paragraph{Results:}  \Cref{2_D_Oscilatory_noise_level} demonstrates the performance of different algorithms in the presence of noise. We consider additive white Gaussian noise with different standard variances $\sigma=\{0,\ 0.02,\ 0.04,\ 0.08\}$. It shows that as we increase the noise level, the \rksindy\ fails to estimate the coefficients. However, \ineuralsindy\ and \deepymod\ are robust in discovering the underlying equations accurately, even for high noise levels, and both exhibit similar performances. In \Cref{2_D_Oscillators_table} (in the appendix), we also report learned governing equations from data with various noise levels, which again illustrate that both \ineuralsindy\ and \deepymod\ have similar performance, and \rksindy~fails to recover governing equations from highly noisy data.
Furthermore, in \Cref{Coeff_iteration_2_D_Oscillator}, the convergence of the non-zero coefficients for the different methods is shown as the training progresses. It can be seen that \ineuralsindy\ has a faster convergence rate compared to \deepymod\ and \rksindy. 
Next, we discuss the performance of \ineuralsindy~and \deepymod~for \scenea~and \sceneb.

\begin{figure}[!tb]
	\centering
	
	\begin{subfigure}[b]{1\textwidth}
	\centering
	\includegraphics[width=0.99\textwidth]{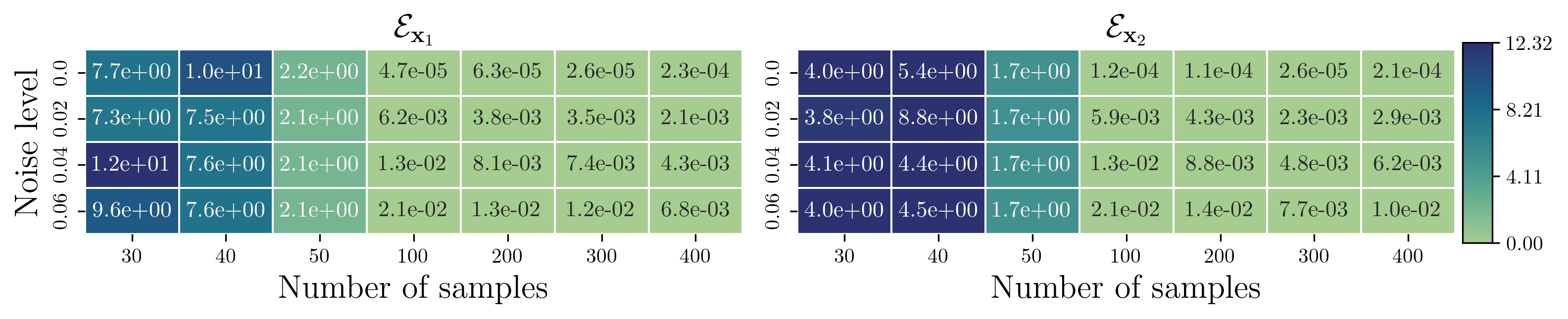}
	\caption{Using \ineuralsindy.}
\end{subfigure}

	\begin{subfigure}[b]{1\textwidth}
		\centering
		\includegraphics[width=0.99\textwidth]{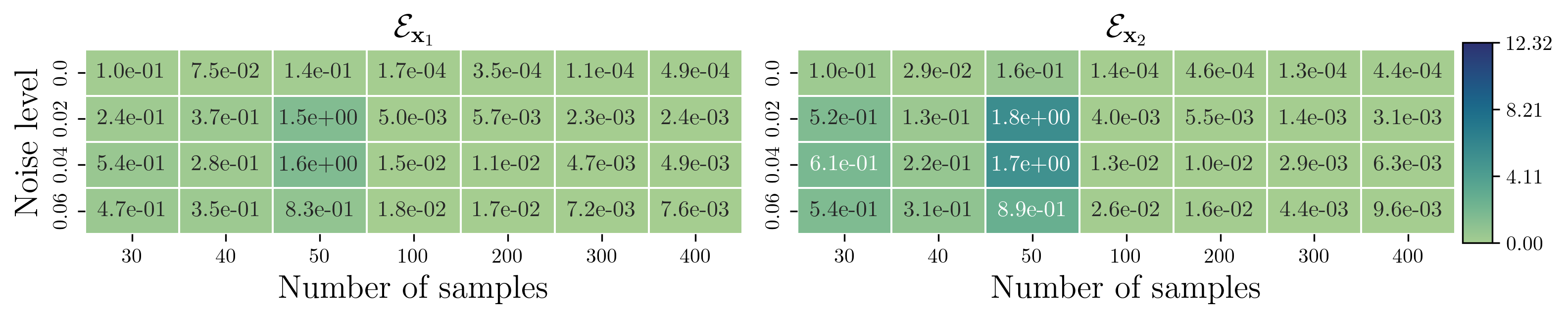}
		\caption{Using \deepymod.}
	\end{subfigure}
	\caption{Linear oscillator: A comparison of \ineuralsindy\ and \deepymod\ under \scenea.}
	\label{heatmap_2_D_Oscillator_fixed_neurons}	
\end{figure}

\begin{figure}[!tb]
	\centering
	
	\begin{subfigure}[b]{1\textwidth}
		\centering
	\includegraphics[width=0.99\textwidth]{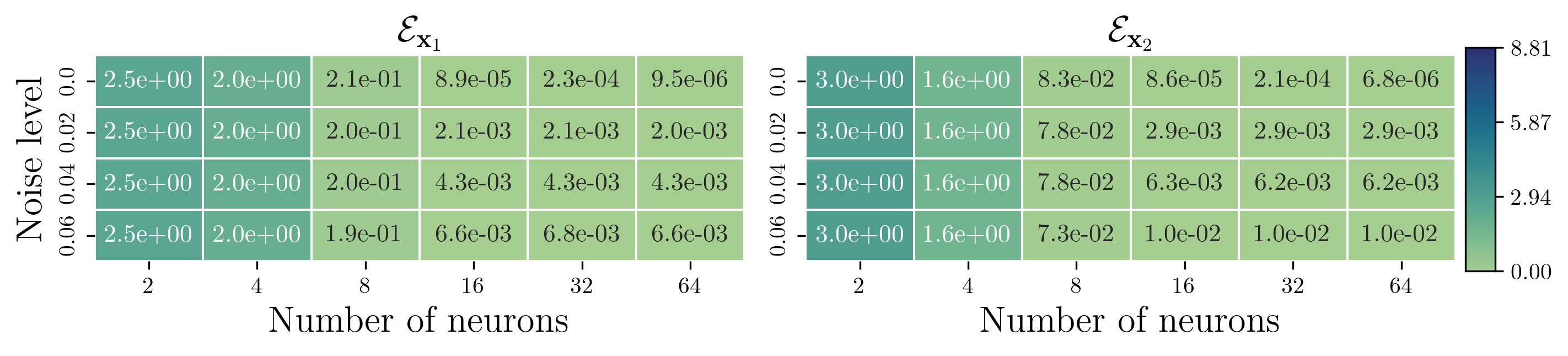}
		\caption{Using \ineuralsindy.}
	\end{subfigure}
	
	\begin{subfigure}[b]{1\textwidth}
		\centering
	\includegraphics[width=0.99\textwidth]{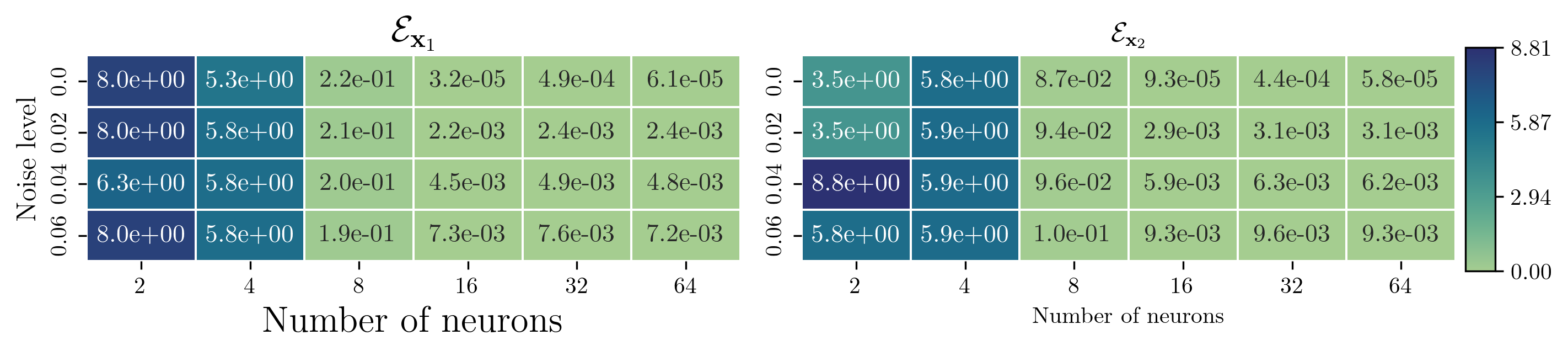}
		\caption{Using \deepymod.}
	\end{subfigure}
	\caption{Linear oscillator: A comparison of \ineuralsindy\ and \deepymod\ under \sceneb}
	\label{heatmap_2_D_Oscillator_fixed_sample}	
\end{figure}

\begin{itemize}
	\item \scenea: We consider a \DNN\ architecture with three hidden layers, each having $32$ neurons. For comparison, we consider noise levels with standard variance $\sigma=\{0,\ 0.02,\ 0.04,\ 0.06\}$, and take the number of samples $\{30,\ 40,\ 50,\ 100,\ 200,\ 300,\ 400\}$ in the time interval $[0,10]$ for a single initial condition $(\bx_1(0), \bx_2(0))=(5, 2)$. The rest of the settings are the same as mentioned earlier in the simulation setup. By varying the noise levels and the number of samples, we report the quality of the learned governing equations in \Cref{heatmap_2_D_Oscillator_fixed_neurons}.  
	Note that the error criterion defined in \eqref{error_criteria} is used. Each cell shows the error corresponding to sample sizes and noise levels. By comparing the simulation results, we notice that \deepymod\ performs better for low data regime, but as the number of data is increased, both \ineuralsindy\ and \deepymod\ perform similarly.
	
	\item \sceneb: In this case,  we consider a \DNN\ architecture with three hidden layers but vary the number of neurons at each layer from $2$ to $64$. Again, we consider various noise levels. We take $400$ samples in the time interval $[0,10]$ for a single arbitrary initial condition $(\bx_1(0),\bx_2(0))=(5, 2)$. The rest of the settings are the same as mentioned earlier in the simulation setup.
	By varying the noise levels and number of neurons, we report a comparison between \ineuralsindy~and \deepymod\ in \Cref{heatmap_2_D_Oscillator_fixed_sample}, where each cell shows the error, corresponding to a specific number of neurons and noise level. These comparisons again show that both methodologies perform comparably and learn correct coefficients with similar performance for a large number of neurons as the \DNN~has more capacities to capture the dynamics present in the data. More interesting, we would like to highlight that both methods do not over-fit as the capacity of the \DNN~is increased. 
\end{itemize}

\subsection{Cubic damped oscillator}
The cubic oscillatory system is given by the following equation:

\begin{equation}
\begin{aligned}
\dot{\bx}_1(t) &= -0.1 {\bx_1^3(t)} + 2.0 {\bx_2^3(t)},\\
\dot{\bx}_2(t) &= -2.0 {\bx_1^3(t)} - 0.1 {\bx_2^3(t)}.
\end{aligned}
\end{equation}
The system consists of two coupled, non-linear differential equations describing the time evolution of two variables, $\bx_1$ and $\bx_2$. Given noisy data, we aim to recover the governing equations and perform a similar analysis as done for the previous example.

\paragraph{Simulation setup:} To generate the training data set, we consider two initial conditions $(\bx_1(0),\bx_2(0))=\{(2, 2), (-2, -2)\}$ and collect $800$ points in the time interval $t \in [0, 10]$. Our \DNN\ architecture has three hidden layers, each having $32$ neurons. We set the number of epoch ${\texttt{max-iter}} = 30,000$, threshold value $\texttt{tol}= 0.05$. The initial training iteration (${\texttt{init-iter}}$) is set to $15,000$ with the learning rate $10^{-4}$ for the \DNN\ parameters and $10^{-3}$ for the coefficient matrix $\Xi^{\texttt{est}}$. After the initial training, for every $q=5,000$ iterations afterward, we employ the sequential thresholding and update the  \DNN\ parameters and $\Xi^{\texttt{est}}$. The dynamical system is estimated in the space of polynomials up to order three.

\begin{figure}[!tb]
	\centering
	\includegraphics[width=0.8\textwidth,trim = 0cm 7cm 0cm 0cm, clip]{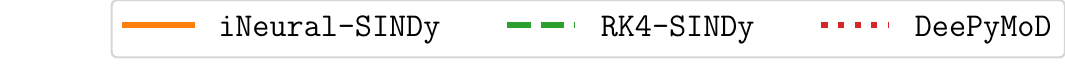}
	\begin{subfigure}{0.49\textwidth}
		\includegraphics[width=1\textwidth,trim = 0cm 0cm 0cm 1.1cm, clip]{\linearexample/Cubic_damped_oscilator/Cubic_Oscilator_noise_0.0_Coefficients__threshold_0.05coeff_1.pdf}
		\caption{Coefficients for $x_1$.}
		\label{}
	\end{subfigure}
	\begin{subfigure}{0.49\textwidth}
		\includegraphics[width=1\textwidth,trim = 0cm 0cm 0cm 1.1cm, clip]{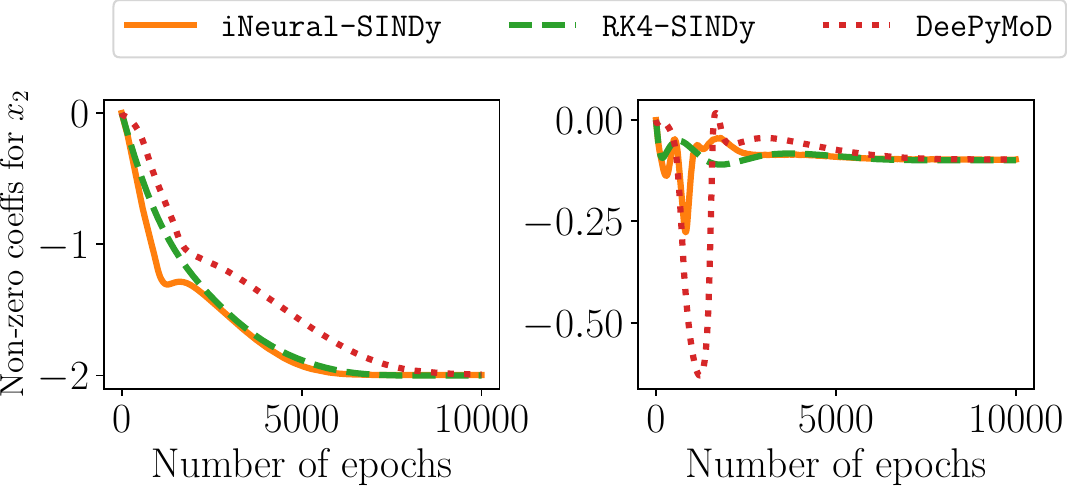}
		\caption{Coefficients for $x_2$.}
		\label{}
	\end{subfigure}
	\caption{Cubic oscillator: Estimated coefficients during the training loop for \ineuralsindy, \deepymod\ and \rksindy.}
	\label{coeff_iteration_cubic_oscillator}
\end{figure}

\begin{figure}[!tb]
	\centering
	
	\begin{subfigure}[b]{1\textwidth}
		\centering
		\includegraphics[width=0.99\textwidth]{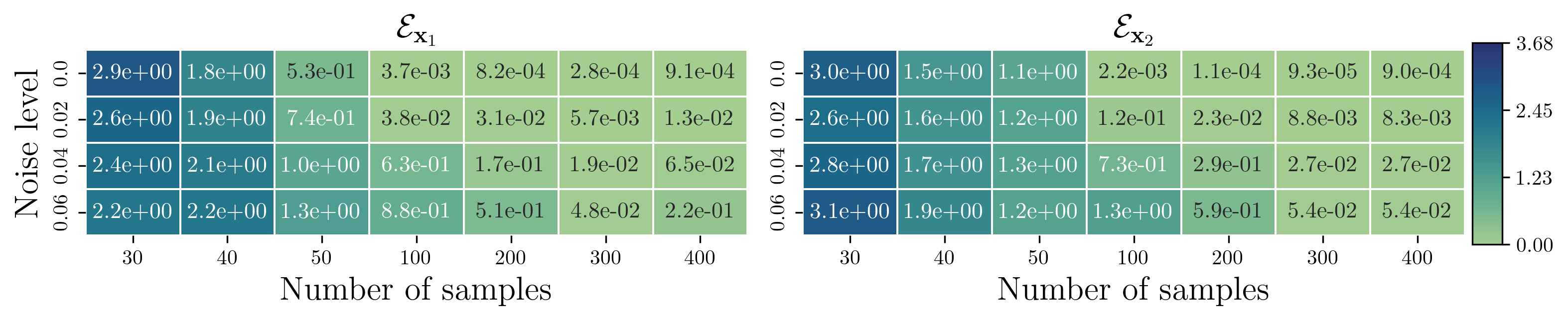}
		\caption{Using \ineuralsindy.}
	\end{subfigure}
	
	\begin{subfigure}[b]{1\textwidth}
		\centering
		\includegraphics[width=0.99\textwidth]{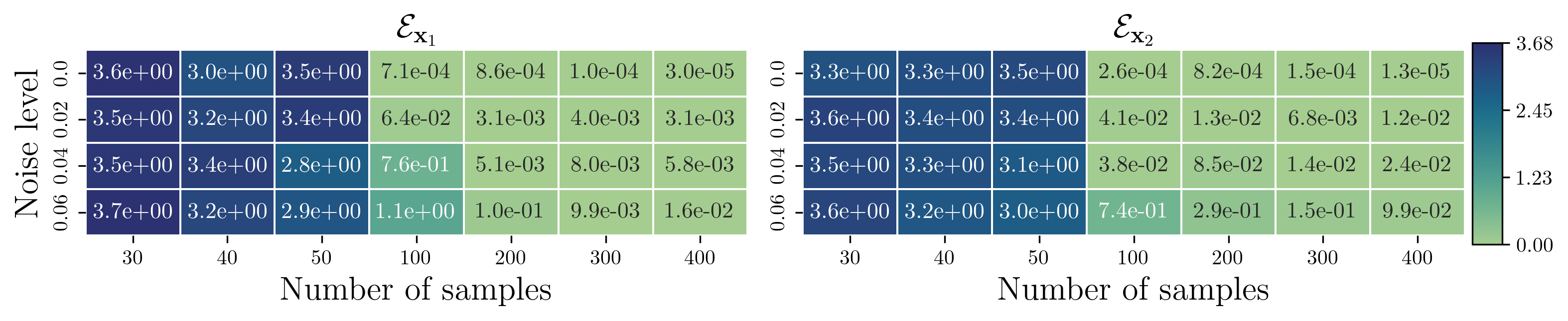}
		\caption{Using \deepymod.}
	\end{subfigure}
	\caption{Cubic oscillator: A comparison of \ineuralsindy\ and \deepymod\ under \scenea.}
	\label{heatmap_Cubic_Oscillator_fixed_neurons}	
\end{figure}

\begin{figure}[!tb]
	\centering
	
	\begin{subfigure}[b]{1\textwidth}
		\centering
		\includegraphics[width=0.99\textwidth]{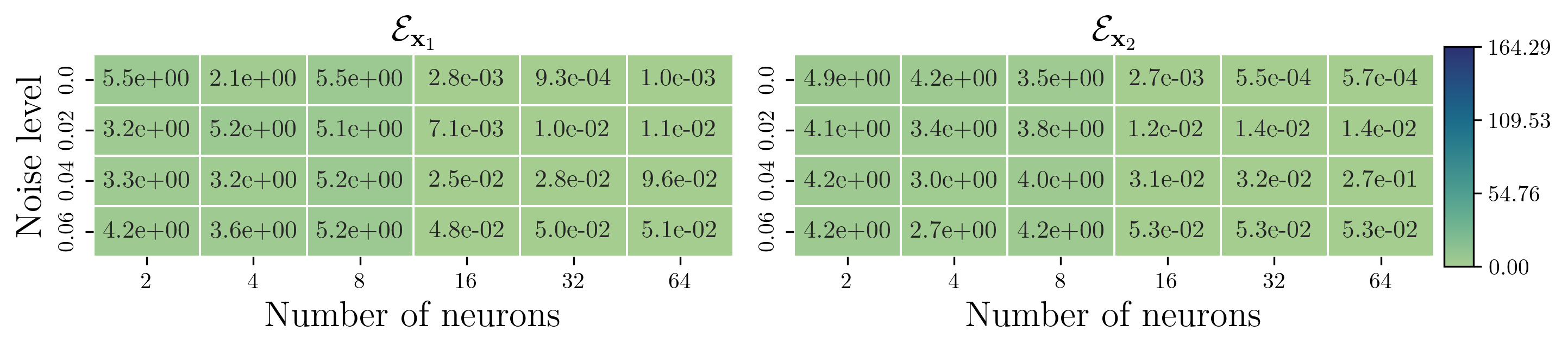}
		\caption{Using \ineuralsindy.}
	\end{subfigure}

	\begin{subfigure}[b]{1\textwidth}
		\centering
		\includegraphics[width=0.99\textwidth]{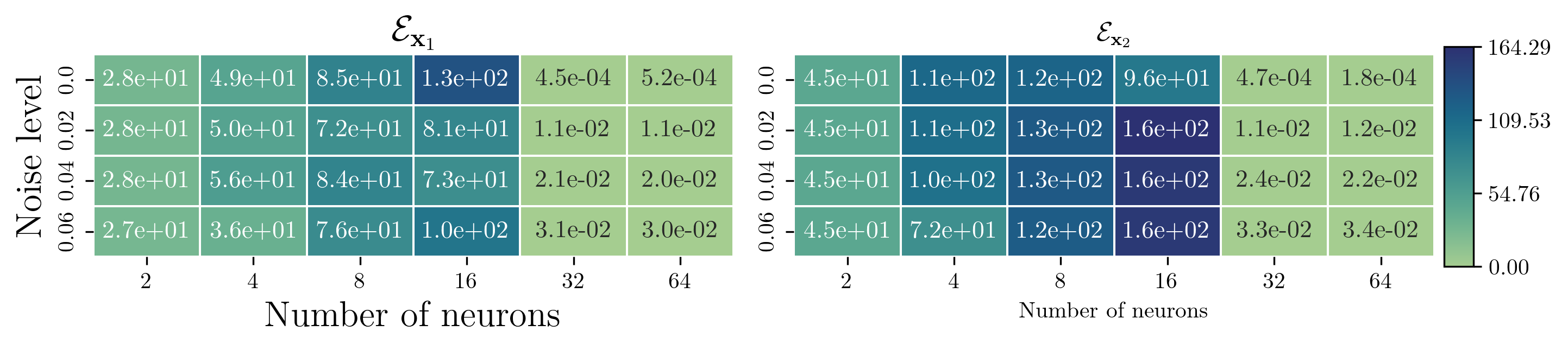}
		\caption{Using \deepymod.}
	\end{subfigure}
	\caption{Cubic oscillator: A comparison of \ineuralsindy\ and \deepymod\ under \sceneb}
	\label{heatmap_Cubic_Oscillator_fixed_samples}	
\end{figure}

\paragraph{Results:}
To see the performance of these different methodologies under the presence of noise, we consider a Gaussian noise with the standard variance $\sigma=\{0,\ 0.02,\ 0.04,\ 0.06\}$. We report the obtained results in \Cref{Cubic_Oscilatory_noise_level} and in \Cref{Cubic oscillator} (see Appendix), and we notice that \rksindy~performs poorly for high noise levels, but \ineuralsindy~and \deepymod~have competitive performance. Further, in \Cref{coeff_iteration_cubic_oscillator}, we plot the convergence of the non-zero coefficients as the training progresses for the noise-free case. Here, we again observe a faster convergence for \ineuralsindy\ as compared to the other two approaches. Next, we investigate performances of \ineuralsindy\ and \deepymod\ for \scenea\ and \sceneb, which are discussed in the following.
\begin{itemize}
	\item \scenea: We fix a \DNN\ architecture with three  hidden layers, each having $32$ neurons. We consider a set of noise level with $\sigma=\{0,\ 0.02,\ 0.04,\ 0.06\}$ and a set of sample size $\{30, 40,\ 50,\ 100,$ $200,\ 300,\ 400\}$. The data are collected using a random initial condition in the interval $[1, 4]$ for $\{\bx_1,\bx_2\}$.
	The rest of the settings are the same as mentioned earlier in the simulation setup. The results are shown in  \Cref{heatmap_Cubic_Oscillator_fixed_neurons}, where we notice that for a smaller data set, \ineuralsindy~performs slightly better as compared to \deepymod, whereas for larger data set, it is otherwise.
	
	\item \sceneb:	For this case, we fix the sample size to $400$ but consider a \DNN\ architecture with three hidden layers with the number of neurons ranging from $2$ to $64$. Furthermore, we consider a set of noise levels with $\sigma=\{0,\ 0.02,\ 0.04,\ 0.06\}$. The data are generated as in \scenea\, and the training setting is also to be as above. The results are depicted in \Cref{heatmap_Cubic_Oscillator_fixed_samples}, where we observe that \ineuralsindy\ performs better as compared to \deepymod~for fewer neurons,  and as we increase the number of neurons, both methods perform similarly.
\end{itemize}

\begin{figure}[!tb]
	\centering
	\includegraphics[width=0.7\textwidth, trim = 0cm 7.8cm 0cm 0cm, clip]{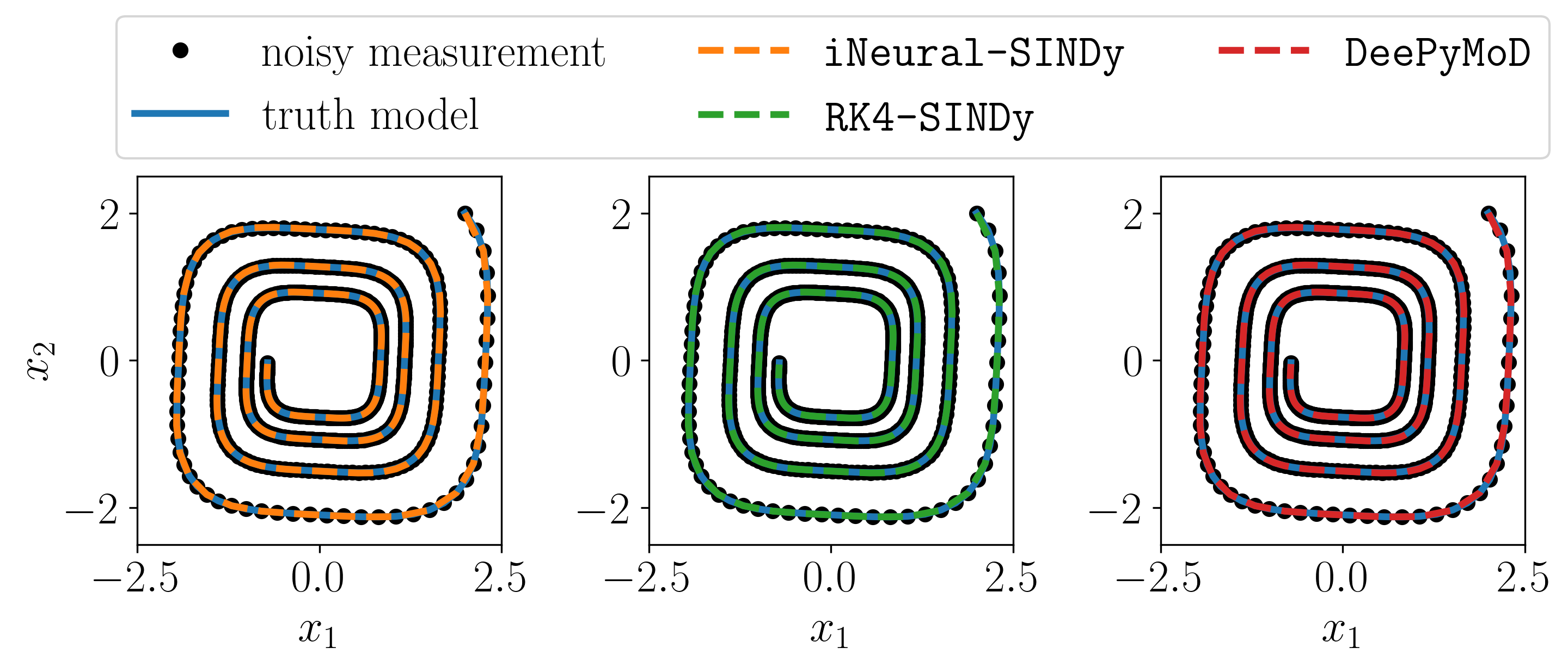} \\[2pt]
	
	\begin{subfigure}{0.49\textwidth}
		\includegraphics[width=\textwidth, trim = 0cm 0cm 0cm 2.75cm, clip]{Pictures/Cubic_Oscilator/Cubic_Oscilator_noise_0.0801_threshold_0.05.png}

		\caption{noise level $\sigma=0.00$}
	\end{subfigure}
	\hfill
	\begin{subfigure}{0.49\textwidth}
		\includegraphics[width=1\textwidth, trim = 0cm 0cm 0cm 2.75cm, clip]{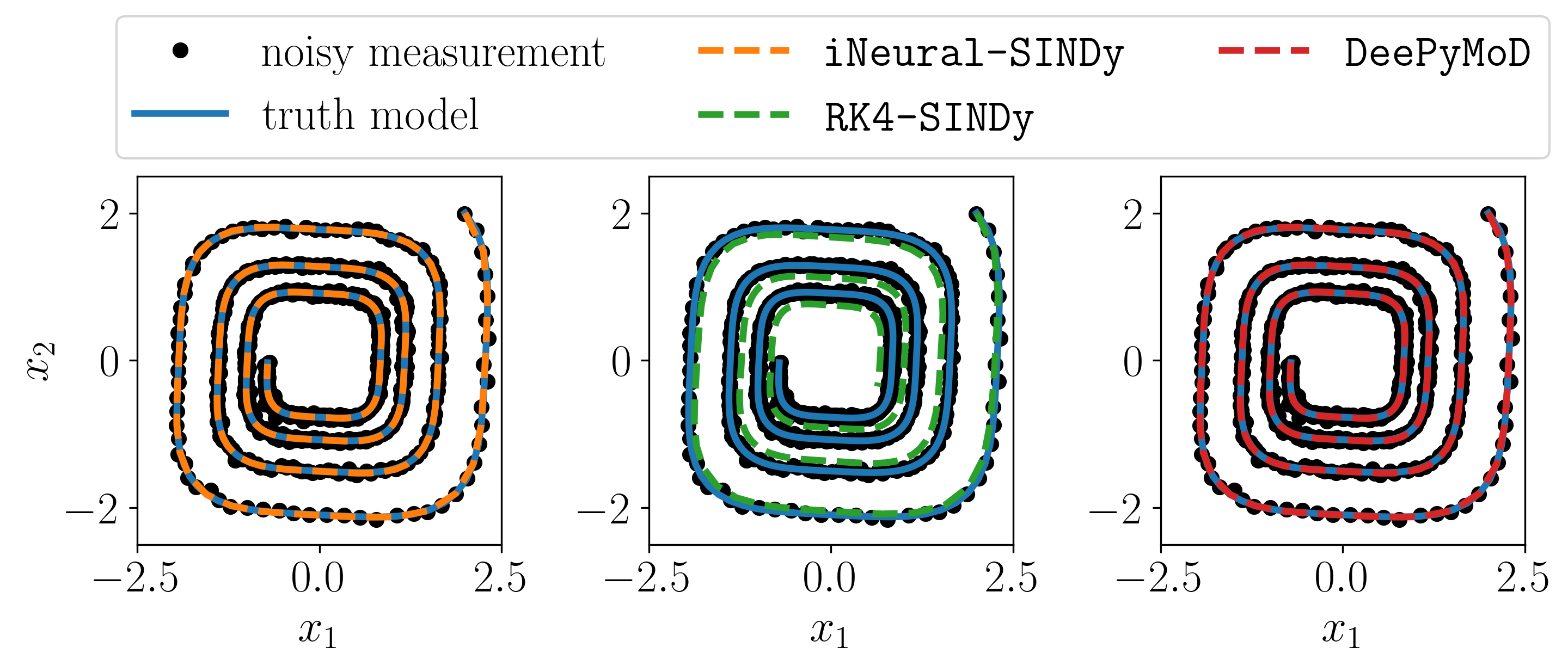}

		\caption{noise level $\sigma=0.02$}
	\end{subfigure}
	\hfill
	\begin{subfigure}{0.49\textwidth}
		\includegraphics[width=1\textwidth, trim = 0cm 0cm 0cm 2.75cm, clip]{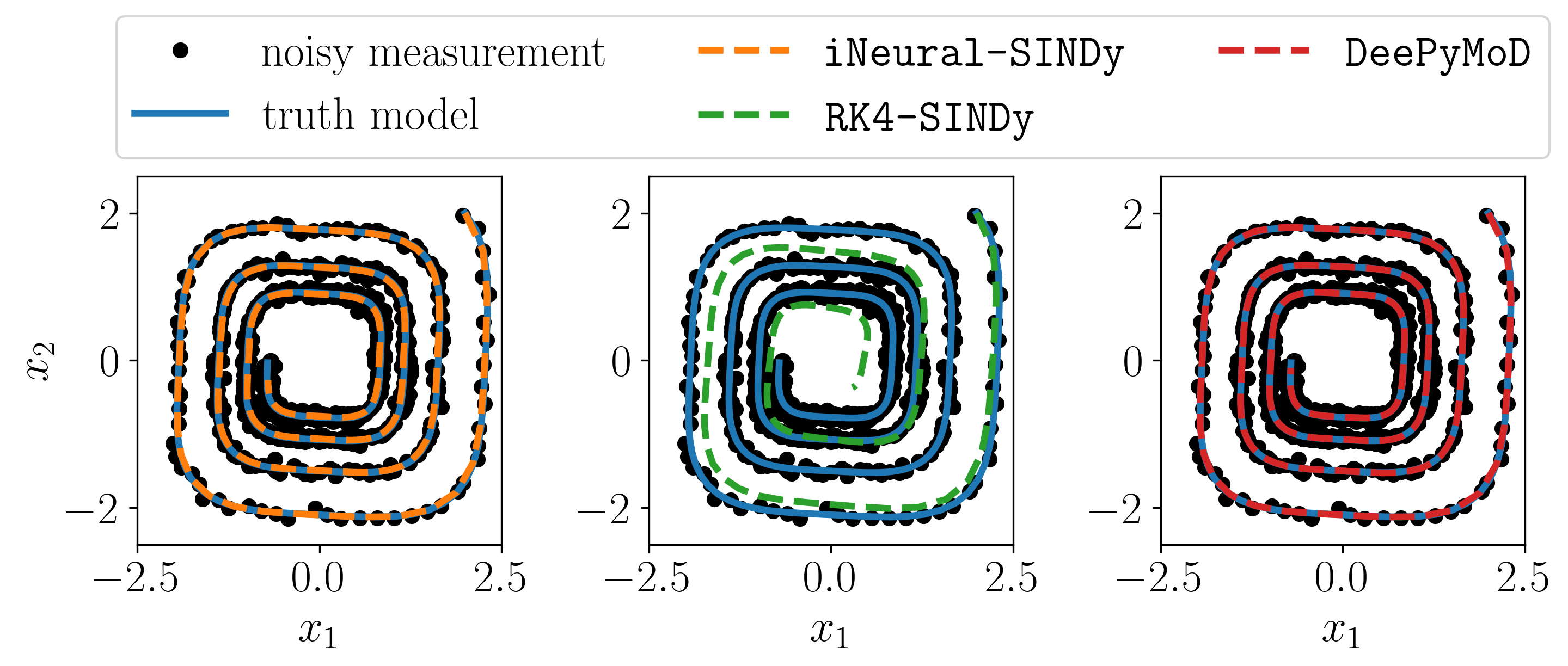}

		\caption{noise level $\sigma= 0.04$}
	\end{subfigure}
	\hfill
	\begin{subfigure}{0.49\textwidth}
		\includegraphics[width=\textwidth,trim = 0cm 0cm 0cm 2.75cm, clip]{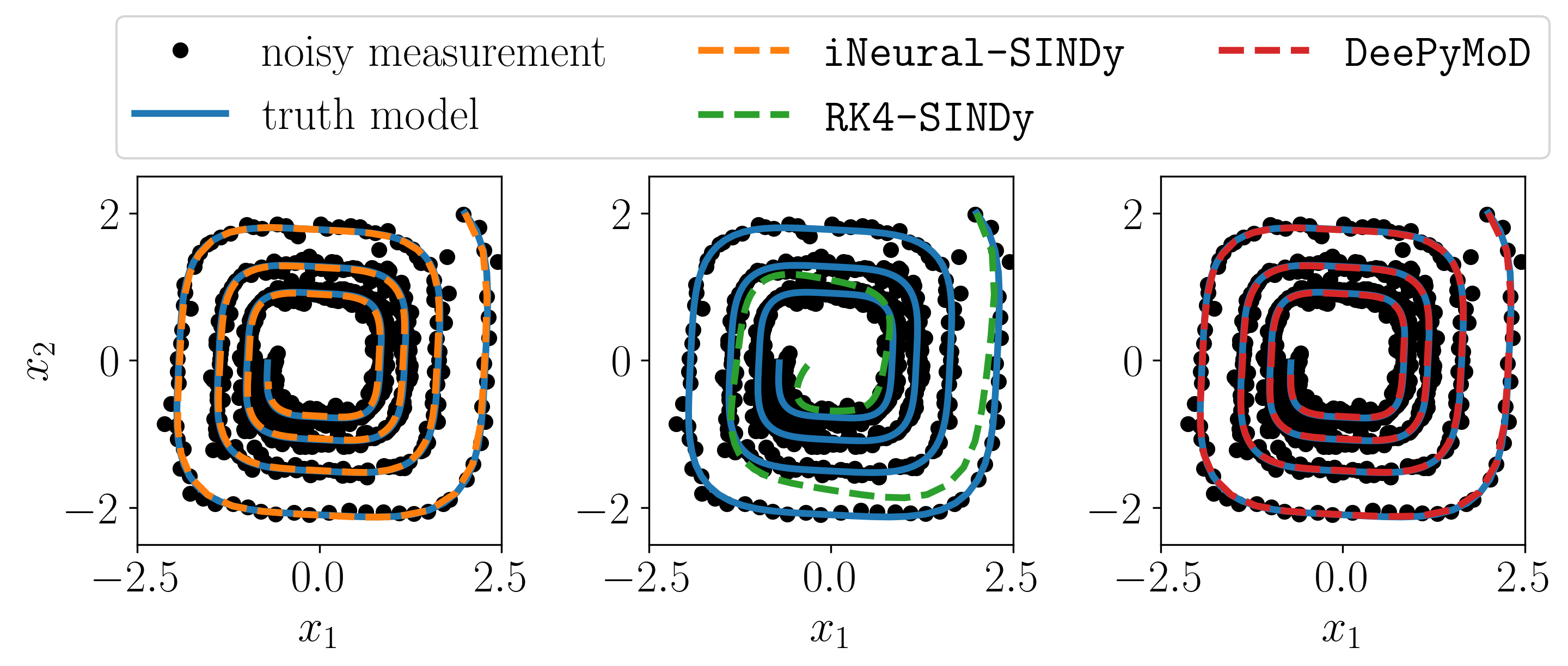}

		\caption{noise level $\sigma = 0.06$}
		\label{Cubic_Oscillatory_noise_06}
	\end{subfigure}
	\caption{Cubic damped oscillator: A comparison of the estimation with different techniques and noise level}
	\label{Cubic_Oscilatory_noise_level}
\end{figure}

\subsection{Fitz-Hugh Nagumo system}
The Fitz-Hugh Nagumo (FHN) model is a non-linear system of ordinary differential equations that is used to describe the behavior of a biological neuron, see, e.g., \cite{Schaeffer_2017}. The FHN model is commonly used to study the behavior of biological neurons under different conditions, such as the changes in the external stimulus or variations in the intrinsic properties of the neuron. The set of differential equations that describe the underlying dynamics are as follows:

\begin{equation}
\begin{aligned}
\dot{\bx}_1(t) & = 1.0{\bx_1(t)} - 1.0 {\bx_2(t)} - \frac{1}{3} {\bx_1^3(t)} + 0.1, \\
\dot{\bx}_2(t) &= 0.1 {\bx_1(t)} - 0.1 {\bx_2(t)}.
\end{aligned}
\end{equation}

\paragraph{Simulation setup:}
For this simulation example, we consider two initial conditions $(\bx_1(0),\bx_2(0))= \{(2, 1.5), (1.5, 2)\} $ and take $400$ data points in the time interval $t \in [0,\ 200]$. The \DNN\ architecture has three hidden layers with $32$ neurons. We set the number of epoch ${\texttt{max-iter}} = 50,000$ and threshold value $\texttt{tol}= 0.05$. The number of iterations for the initial training is set to $15,000$ with the learning rate $10^{-4}$ for the \DNN\ parameters and $10^{-3}$ for the coefficient matrix $\Xi^{\texttt{est}}$. After the initial training, we employ the sequential thresholding after each $q=5,000$ iterations. We aim to learn the underlying governing equations in the space of polynomials with degrees up to order three. 

\begin{figure}[!tb]
	\centering
	\includegraphics[width=0.9\textwidth,trim = 0cm 5.5cm 0cm 0cm, clip]{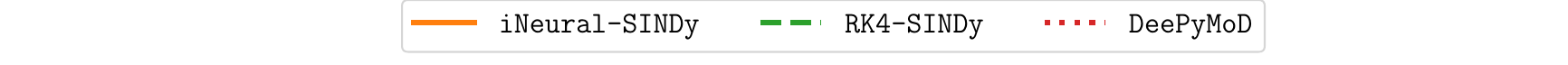}\\[2pt]
	
	\begin{subfigure}{1\textwidth}
		\includegraphics[width=1\textwidth,trim = 0cm 0cm 0cm 1.1cm, clip]{\linearexample/Fitz_hugh_nagumo/Fitz_Hugh_Nagumo_noise_0.0_Coefficients__threshold_0.05coeff_1.pdf}
		\caption{Coefficients for $x_1$.}
		\label{}
	\end{subfigure}
	\begin{subfigure}{1\textwidth}
		\includegraphics[width=1\textwidth,trim = 0cm 0cm 0cm 1.1cm, clip]{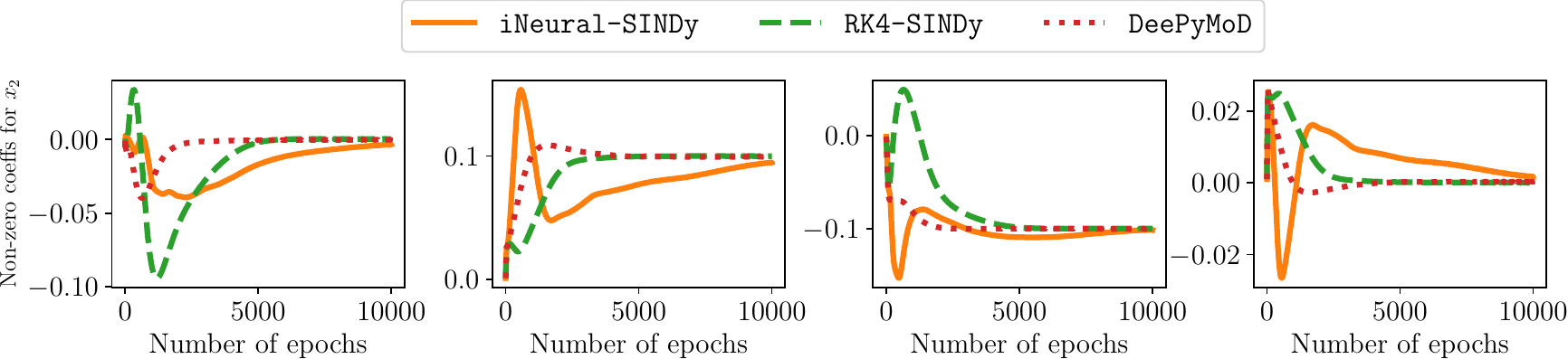}
		\caption{Coefficients for $x_2$.}
		\label{}
	\end{subfigure}
	\caption{Fitz-Hugh Nagumo: Estimated coefficients during the training loop for \ineuralsindy, \deepymod\ and \rksindy.}
	\label{coeff_iteration_Fitz-Hugh Nagumo}
\end{figure}

\paragraph{Results:}
Converse to the results that we earned in the previous two examples, for the FHN, \ineuralsindy\ has a slower convergence rate compared to \deepymod~ and \rksindy, see \Cref{coeff_iteration_Fitz-Hugh Nagumo}.

For this example, we again make a similar observation (see \Cref{Fitz_Hugh_Nagumo_noise_level}, and \Cref{Fitz-Hugh_Nagumo_table} in Appendix), where we notice that \ineuralsindy~and \deepymod~exhibit similar performances for lower noise levels, but for the higher noise values, (see the results for $\sigma = 0.08$ \Cref{Fitz-Hugh_Nagumo_table}), \ineuralsindy~tends to outperform \deepymod. Moreover, \rksindy~clearly fails for high noise levels. However, converse to the results reported in the previous two examples, for this example, we notice a slower convergence of \ineuralsindy\ compared to \deepymod~ and \rksindy; see, \Cref{coeff_iteration_Fitz-Hugh Nagumo}. But we highlight that \ineuralsindy~can identify governing equations for highly noisy data, as stated earlier.
Next, we compare the performances of \ineuralsindy~and \deepymod~under \scenea~and \sceneb.

 \color{black}
\begin{figure}[!tb]
	\centering
	\includegraphics[width=0.7\textwidth, trim = 0cm 7.8cm 0cm 0cm, clip]{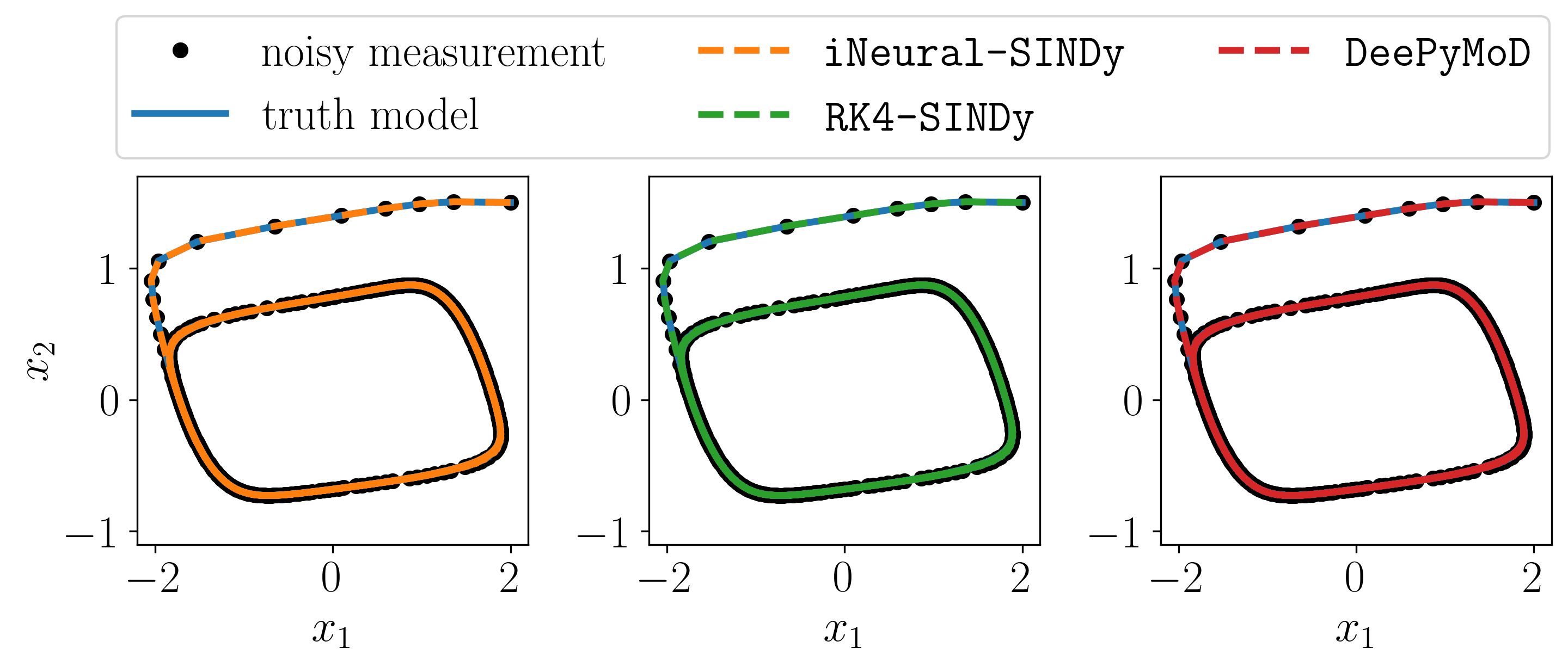} \\[2pt]

	\begin{subfigure}{0.49\textwidth}
		\includegraphics[width=\textwidth,trim = 0cm 0cm 0cm 2.75cm, clip]{Pictures/Fitz_Hugh_Nagumo/Fitz_Hugh_Nagumo_noise_0.0401_threshold_0.05.png}
		\caption{noise level $\sigma= 0.00$}
	\end{subfigure}
	\hfill
	\begin{subfigure}{0.49\textwidth}
		\includegraphics[width=\textwidth,trim = 0cm 0cm 0cm 2.75cm, clip]{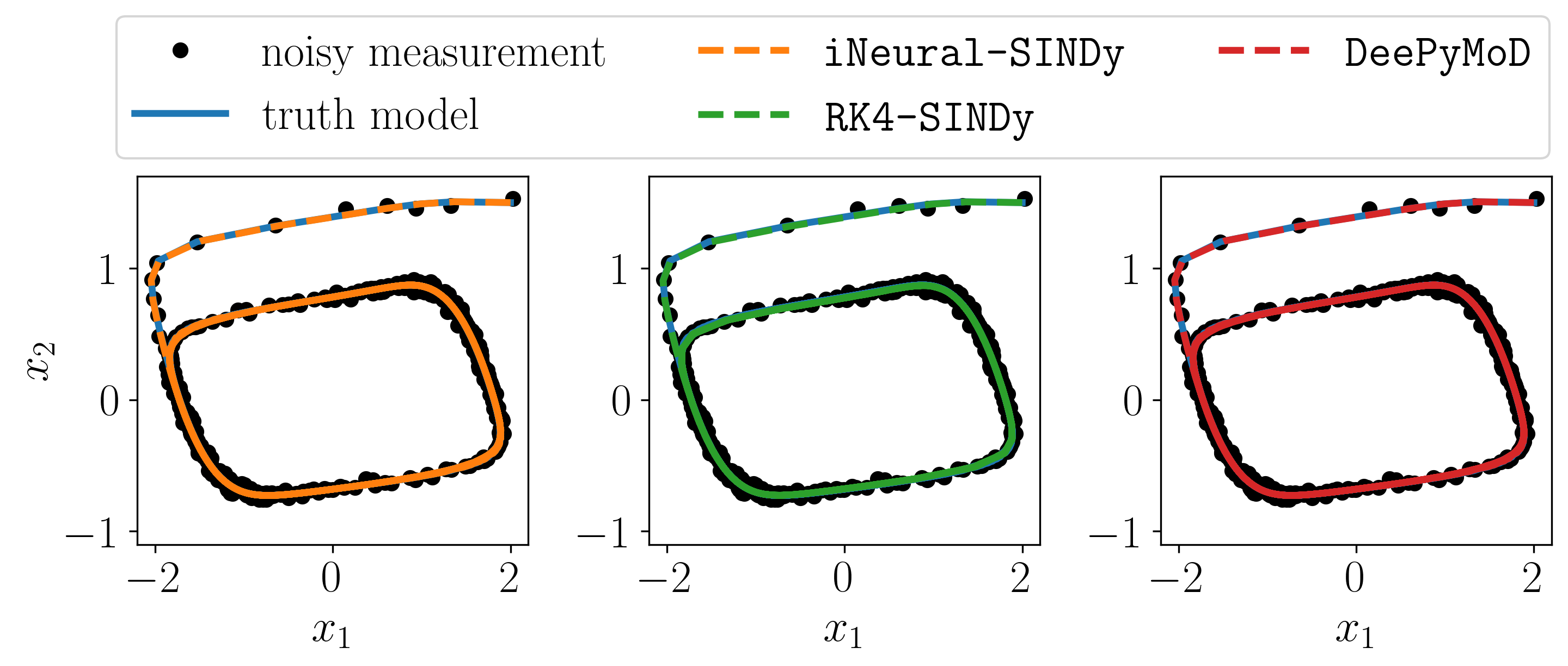}
		\caption{noise level $\sigma= 0.02$}
	\end{subfigure}
	\hfill
	\begin{subfigure}{0.49\textwidth}
		\includegraphics[width=\textwidth,trim = 0cm 0cm 0cm 2.75cm, clip]{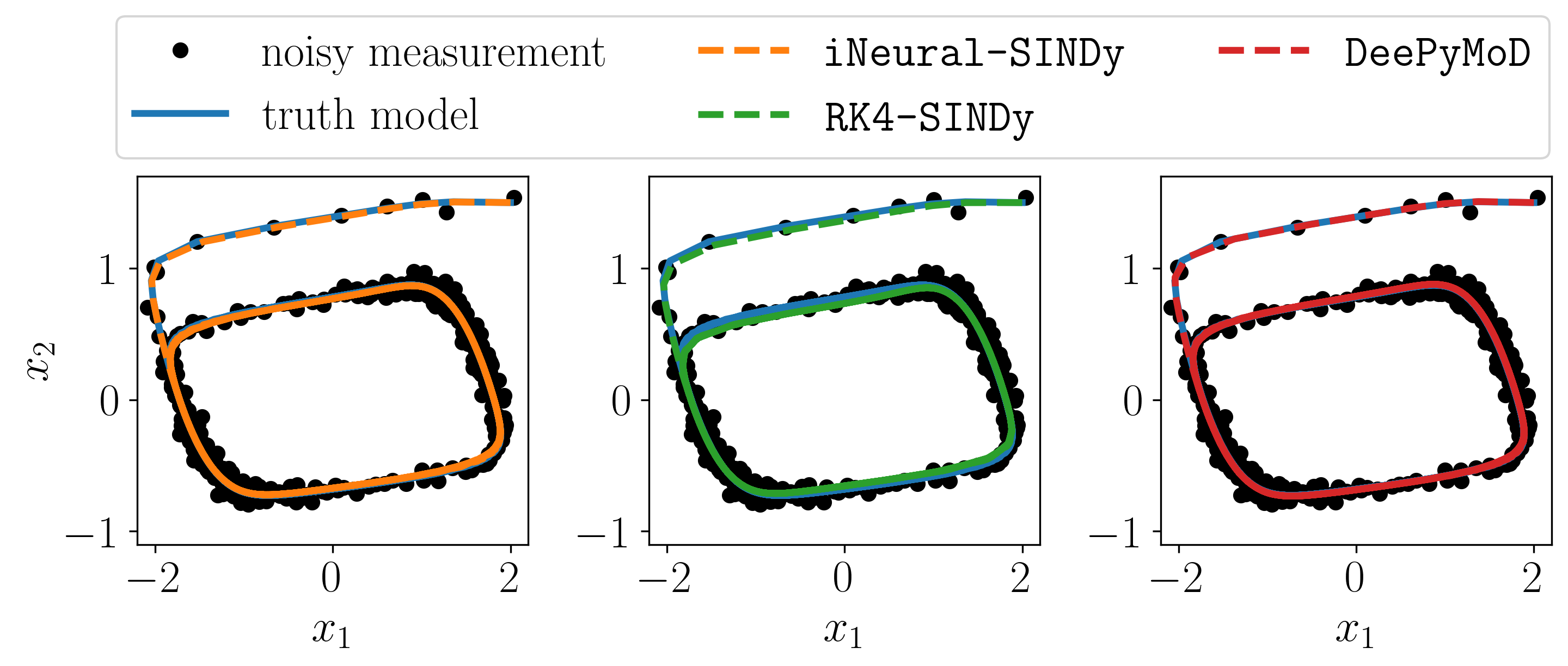}
		\caption{noise level $\sigma= 0.04$}
	\end{subfigure}
	\hfill
	\begin{subfigure}{0.49\textwidth}
		\includegraphics[width=\textwidth,trim = 0cm 0cm 0cm 2.75cm, clip]{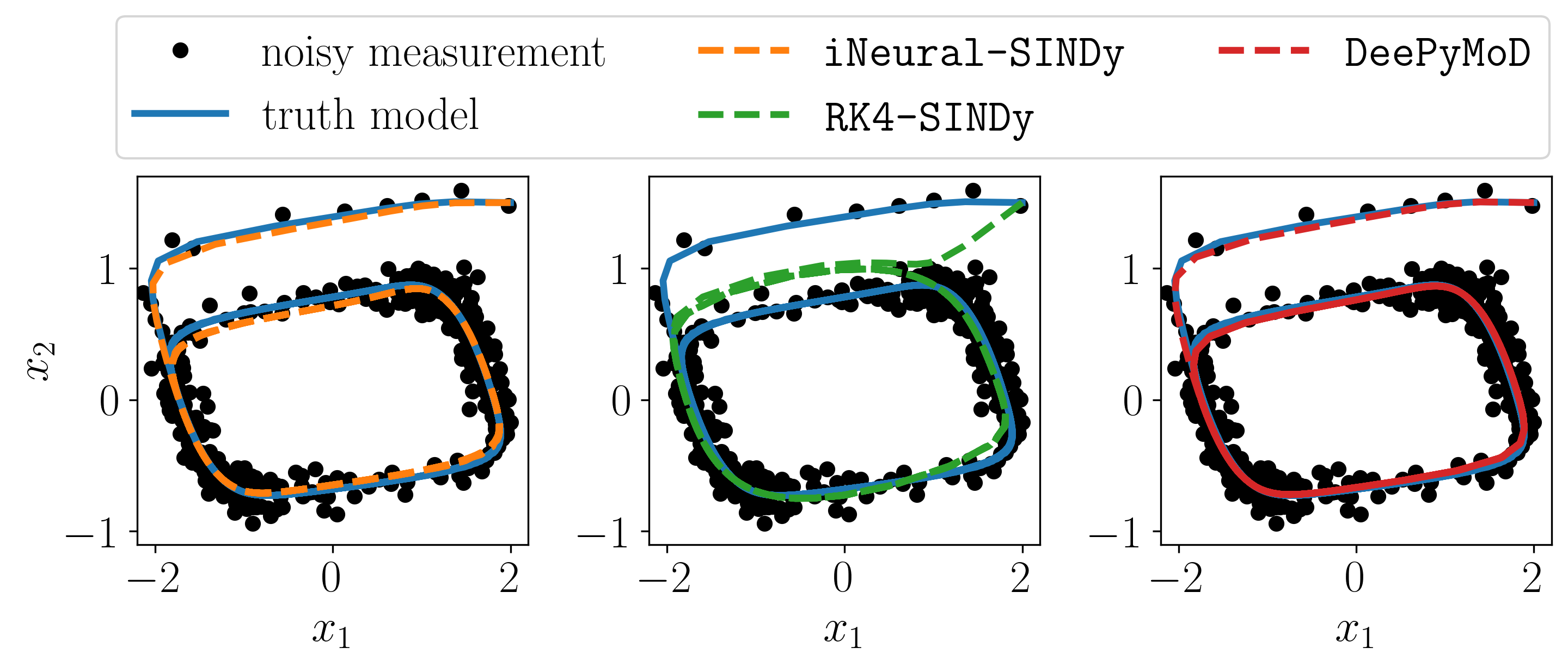}
		\caption{noise level $\sigma= 0.08$}
	\end{subfigure}
	\caption{Fitz-Hugh Nagumo: Comparison of the estimation with different techniques and noise level}
	\label{Fitz_Hugh_Nagumo_noise_level}
\end{figure}

\begin{figure}[!tb]
	\centering
	
	\begin{subfigure}[b]{1\textwidth}
		\centering
		\includegraphics[width=0.99\textwidth]{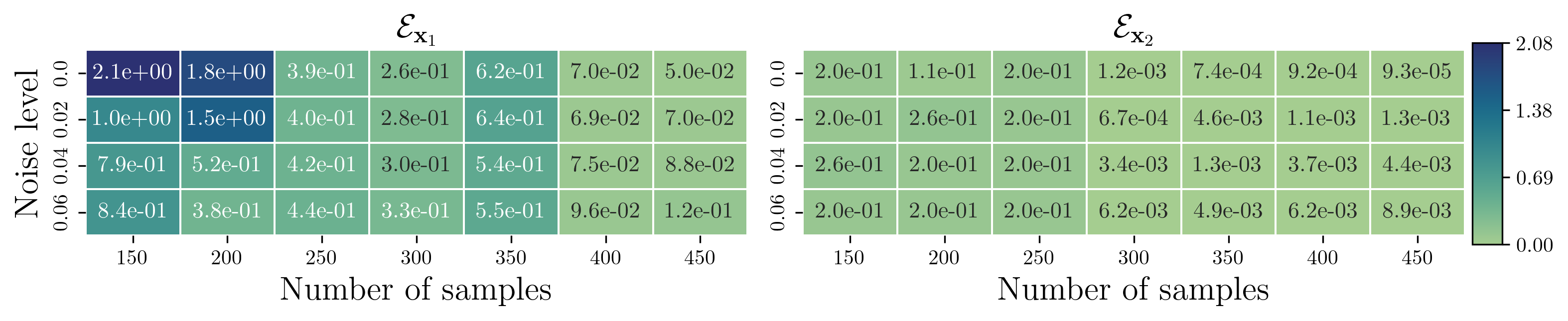}
		\caption{Using \ineuralsindy.}
	\end{subfigure}
	
	\begin{subfigure}[b]{1\textwidth}
		\centering
		\includegraphics[width=0.99\textwidth]{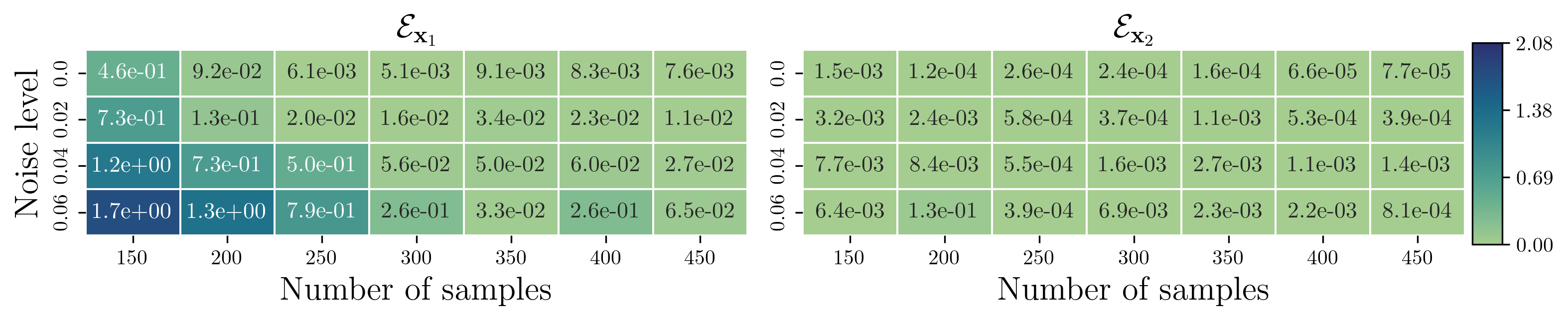}
		\caption{Using \deepymod.}
	\end{subfigure}
	\caption{Fitz-Hugh Nagumo: A comparison of \ineuralsindy\ and \deepymod\ under \scenea.}
	\label{heatmap_Fitz-Hugh Nagumo_fixed_neurons}	
\end{figure}

\begin{figure}[!tb]
	\centering
	
	\begin{subfigure}[b]{1\textwidth}
		\centering
		\includegraphics[width=0.99\textwidth]{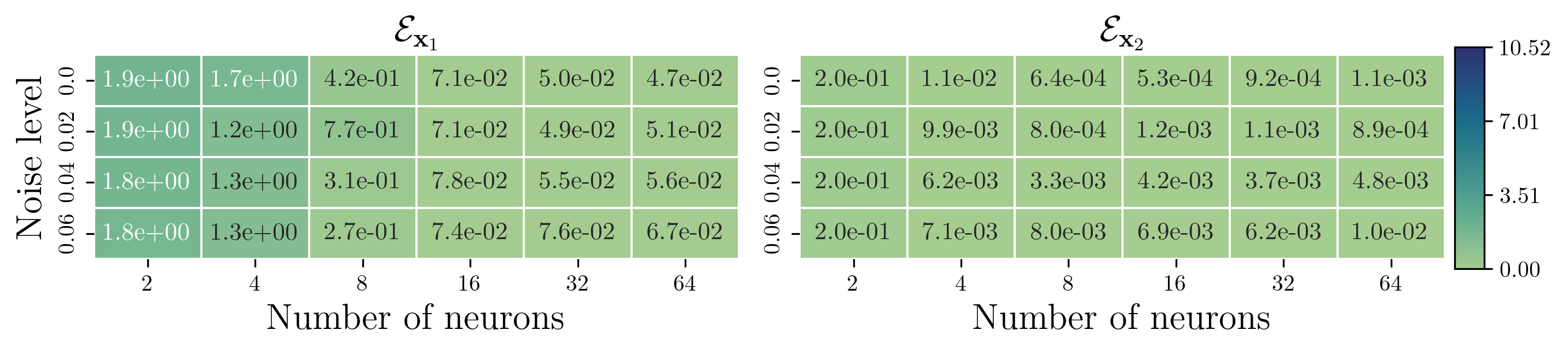}
		\caption{Using \ineuralsindy.}
	\end{subfigure}
	
	\begin{subfigure}[b]{1\textwidth}
		\centering
		\includegraphics[width=0.99\textwidth]{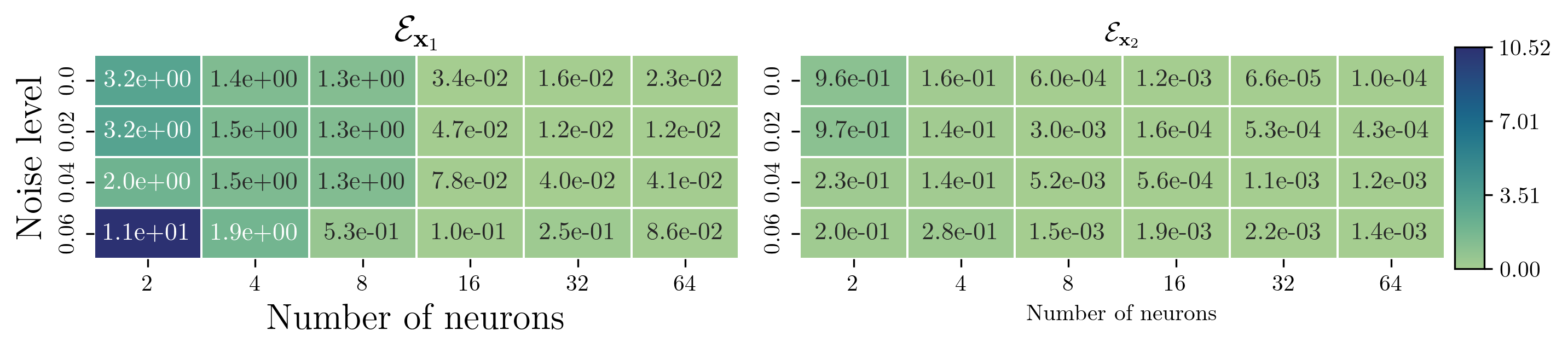}
		\caption{Using \deepymod.}
	\end{subfigure}
	\caption{Fitz-Hugh Nagumo: A comparison of \ineuralsindy\ and \deepymod\ under \sceneb}
	\label{heatmap_Fitz-Hugh Nagumo_fixed_sampls}	
\end{figure}

\begin{itemize}
	\item \scenea: In this case, we consider a fixed \DNN\ architecture with three hidden layers, each consisting of $32$ neurons.  The different noise levels $\{0.0,\ 0.02,\ 0.04, 0.06\}$ are considered, while the sample size is considered in the range from $150$ to $450$ with an increment of $50$. Here, a single initial condition is used for data collection; that is, $(\bx_1(0),\bx_2(0))=(3,2)$. The rest of the settings are the same as mentioned earlier for this example. The results are shown in \Cref{heatmap_Fitz-Hugh Nagumo_fixed_neurons}, where we notice that \deepymod\ outperforms \ineuralsindy\ and has a better performance.
	
	\item 	\sceneb: Here, we conduct a study where we keep the number of samples fixed at $400$, obtained using the initial condition $(\bx_1(0),\bx_2(0))=(3, 2)$. The \DNN\ architecture is designed to have three hidden layers. We aim to explore how \ineuralsindy\ and \deepymod\ perform under different combinations of neurons for each layer and noise level. The training settings for each case remain the same, as mentioned earlier. The outcomes are presented in the heat-map depicted in \Cref{heatmap_Fitz-Hugh Nagumo_fixed_sampls}, where we notice that both \deepymod\ and \ineuralsindy\ almost have the same performance in all the settings.
\end{itemize}

\subsection{Chaotic Lorenz system}\label{subsec:lorenzsetip}
The chaotic Lorenz system is a set of three differential equations as follows \cite{kuznetsov2020attractor}:
\begin{subequations}\label{lorenz_sys}
\begin{align}
\dot{\bx}_1(t) & = \gamma \big( \bx_2(t) - \bx_1(t) \big), \label{lorenz_sys_a}\\
\dot{\bx}_2(t) & = \bx_1(t) \big(\rho - \bx_3(t) \big) - \bx_2(t), \label{lorenz_sys_b}\\
\dot{\bx}_3(t) & = \bx_1(t) \bx_2(t) - \beta \bx_3(t), \label{lorenz_sys_c}
\end{align}
\end{subequations}
where the parameters $\gamma$, $\rho$, and $\beta$ are positive constants with associated standard values $\gamma=10,\ \rho=28,\ \beta=\frac{8}{3}$. The Lorenz system is a classic example of a chaotic system, which means that small differences in the initial conditions can lead to vastly different outcomes over time. It is a widely used benchmark example for discovering governing equations \cite{brunton2016discovering}.

\paragraph{Simulation setup:}
We collect our data in the time interval $t\in [0, 10]$ with a sample size of $200$ for three different initial conditions $(\bx_1(0),\bx_2(0),\bx_3(0))=\{ (-8, 7, 27), (-6, 6, 25), (-9, 8, 22)\}$. The \DNN\ architecture has three hidden layers, each having $64$ neurons. We set the number of iterations ${\texttt{max-iter}} = 35,000$, and threshold value $\texttt{tol}= 0.2$. We set the number of iterations for the initial training to $\texttt{init-iter} = 10,000$, and the learning rate $7 \cdot 10^{-4}$ for the \DNN\ parameters and $10^{-2}$ for the coefficient matrix $\Xi^{\texttt{est}}$. After finishing the initial iterations, we employ sequential thresholding after every $q=3,000$ iterations. Moreover, after each sequential thresholding, we reset the learning rate for \DNN\ parameters to $5 \cdot 10^{-6}$ and for the coefficient matrix $\Xi^{\texttt{est}}$ to $10^{-2}$. The governing equations are estimated by constructing a dictionary with polynomials up to degree two. 

Since the magnitude of $\{\bx_1,\bx_2,\bx_3\}$ for the Lorenz example can be large, we consider scaling the $\bx$'s using a scaling factor $\alpha$. Note that such scaling does not affect the interaction between different $\bx$'s; thus, the sparsity pattern remains the same as well. However, it is observed that improving the condition number of the dictionary matrix enhances the estimate of the coefficients and helps us to determine the right governing equations.

\begin{figure}[!tb]
	\centering
	
	\begin{subfigure}[b]{0.49\textwidth}
		\centering
		\includegraphics[width=1\textwidth]{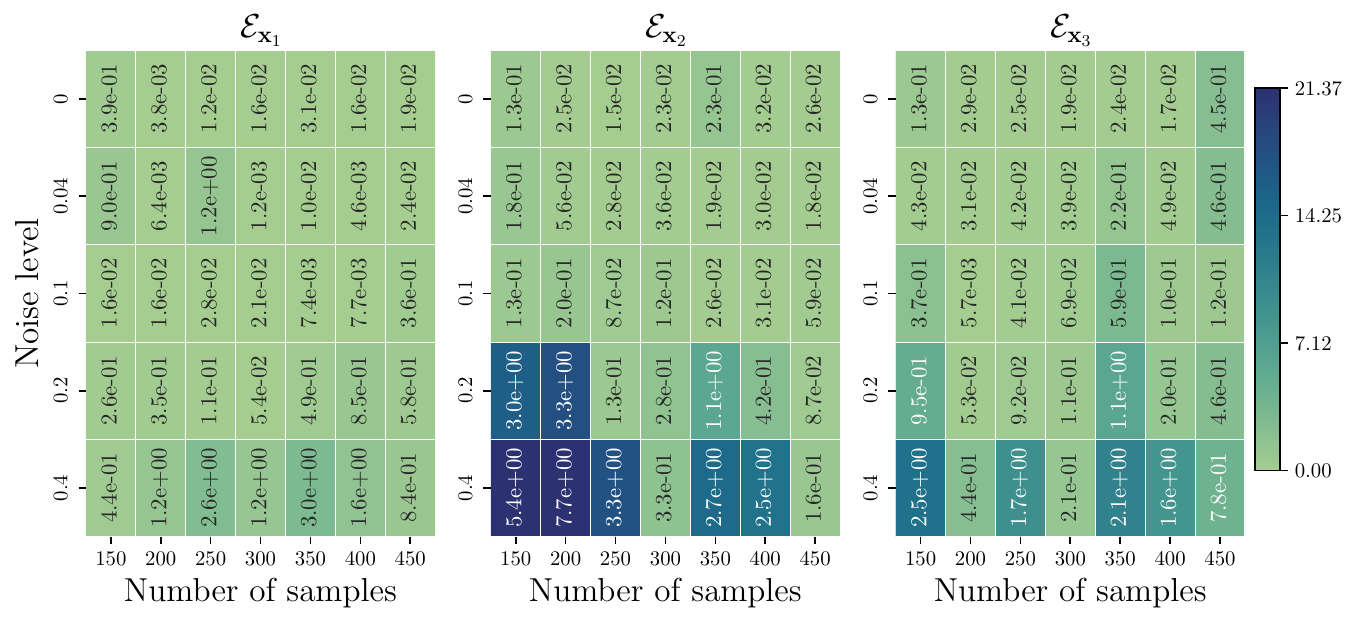}
		\caption{Using \ineuralsindy.}
	\end{subfigure}
	\begin{subfigure}[b]{0.49\textwidth}
		\centering
		\includegraphics[width=1\textwidth]{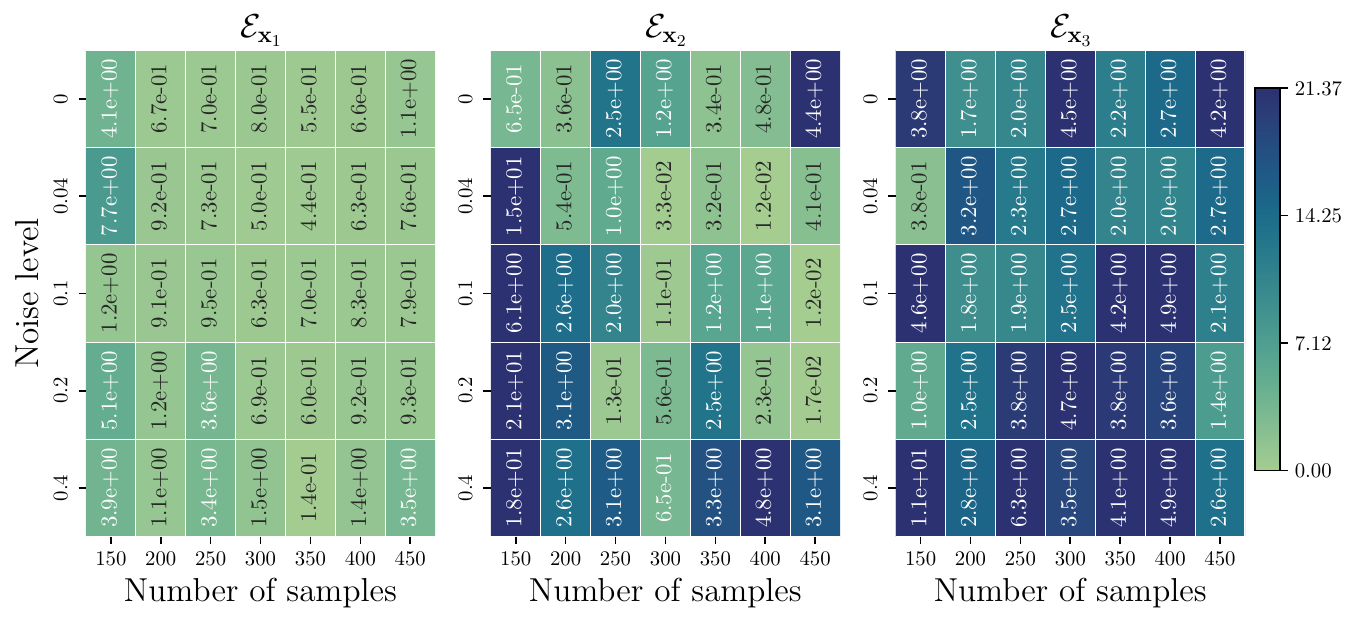}
		\caption{Using \deepymod.}
	\end{subfigure}
	\caption{Lorenz example: A comparison of \ineuralsindy\ and \deepymod\ under \scenea\ with the scaling factor $\alpha=0.1$}
	\label{heatmap_lorenz_fixed_neurons_sacling_0.1}	
\end{figure}

\begin{figure}[!tb]
	\centering
	
	\begin{subfigure}[b]{0.49\textwidth}
		\centering
		\includegraphics[width=1\textwidth]{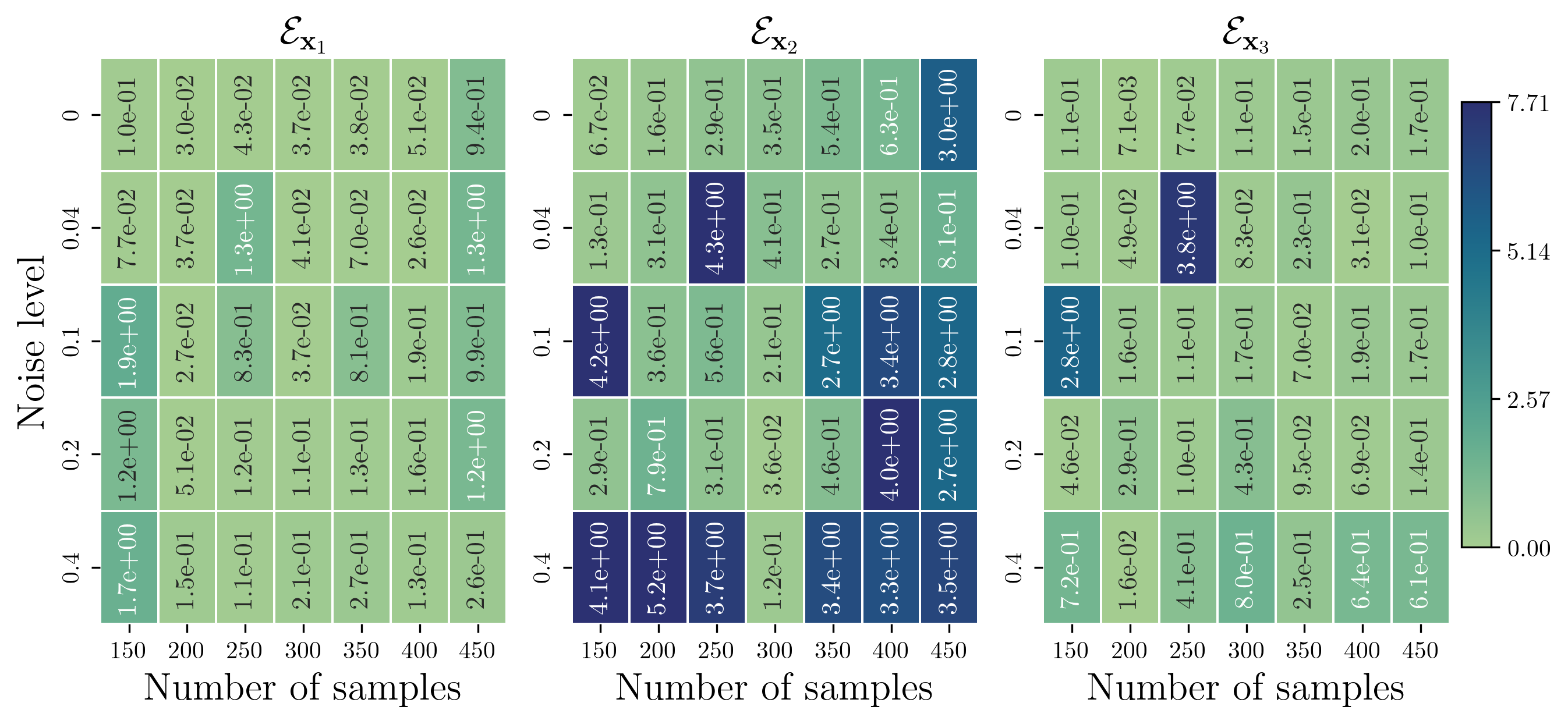}
		\caption{Using \ineuralsindy.}
	\end{subfigure}
	\begin{subfigure}[b]{0.49\textwidth}
		\centering
		\includegraphics[width=1\textwidth]{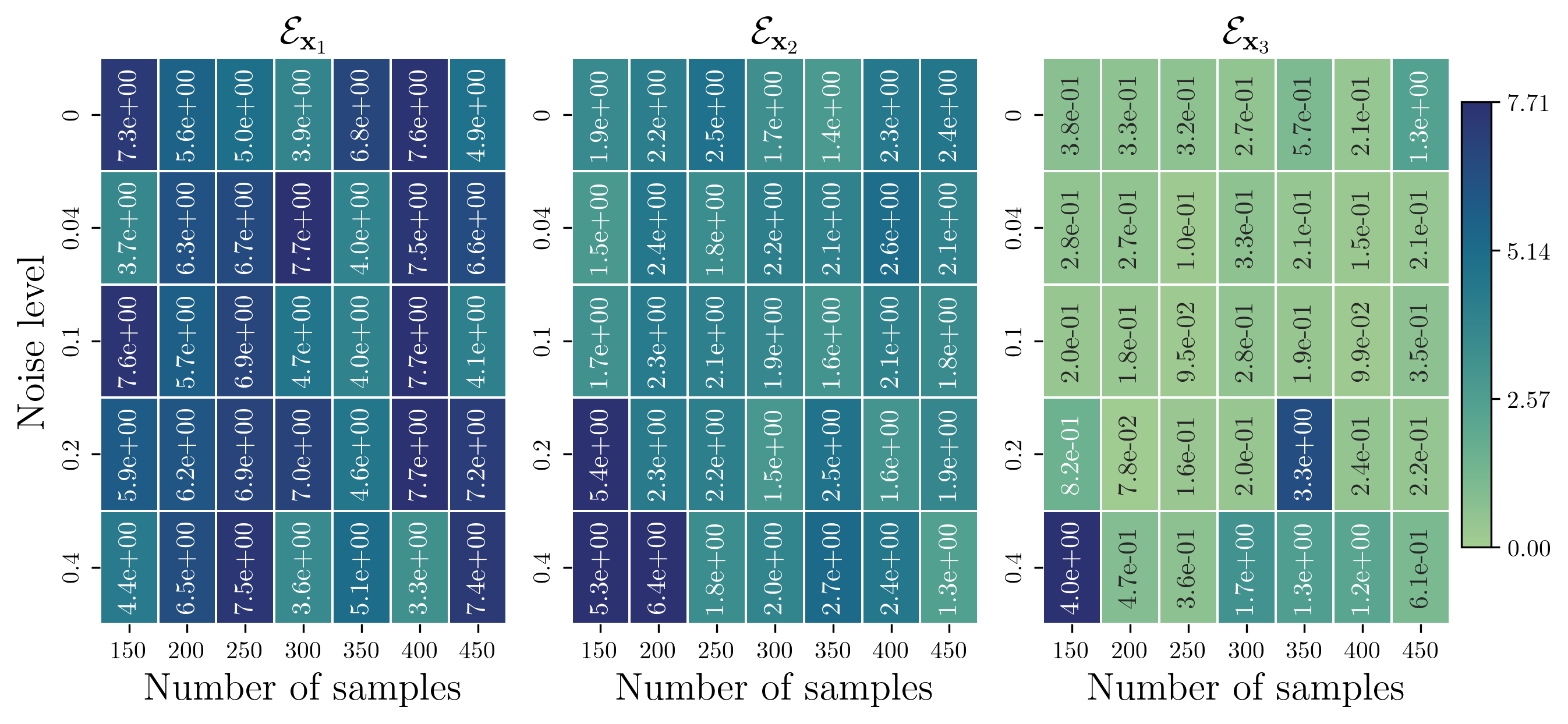}
		\caption{Using \deepymod.}
	\end{subfigure}
	\caption{Lorenz example: A comparison of \ineuralsindy\ and \deepymod\ under \scenea\ with the scaling factor $\alpha=1$}
	\label{heatmap_lorenz_fixed_neurons_sacling_1}	
\end{figure}

\begin{figure}[!tb]
	\centering
	
	\begin{subfigure}[b]{0.49\textwidth}
		\centering
		\includegraphics[width=1\textwidth]{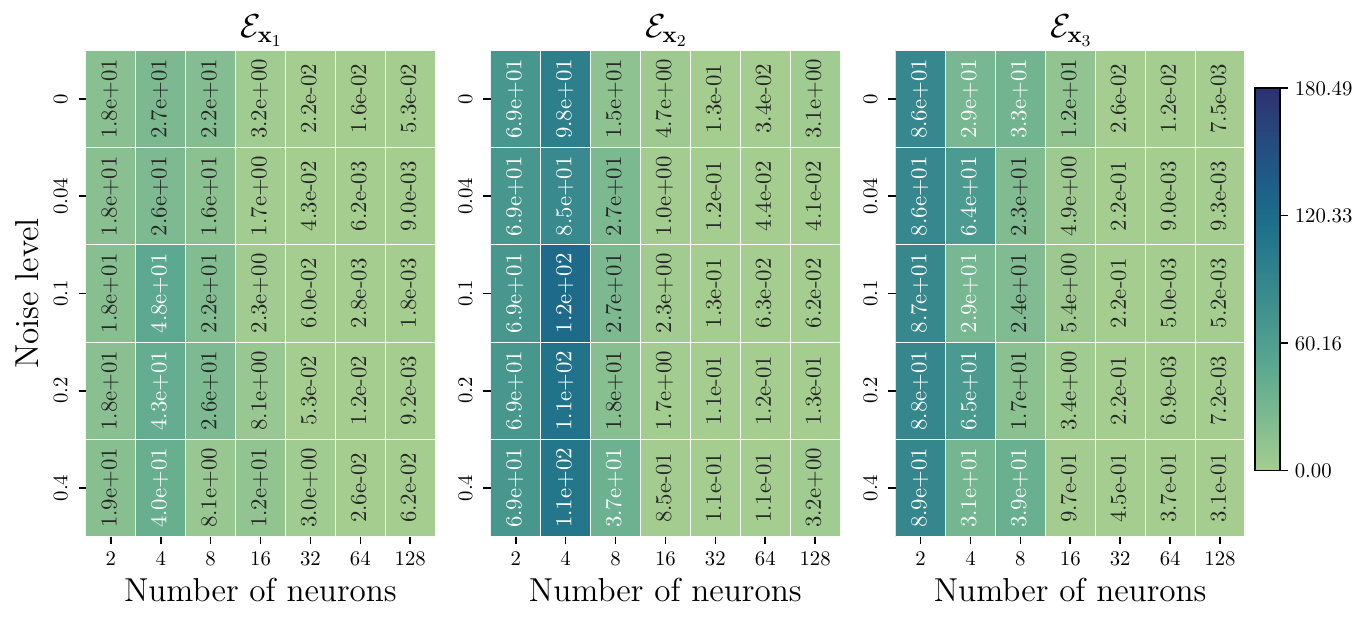}
		\caption{Using \ineuralsindy.}
	\end{subfigure}
	\begin{subfigure}[b]{0.49\textwidth}
		\centering
		\includegraphics[width=1\textwidth]{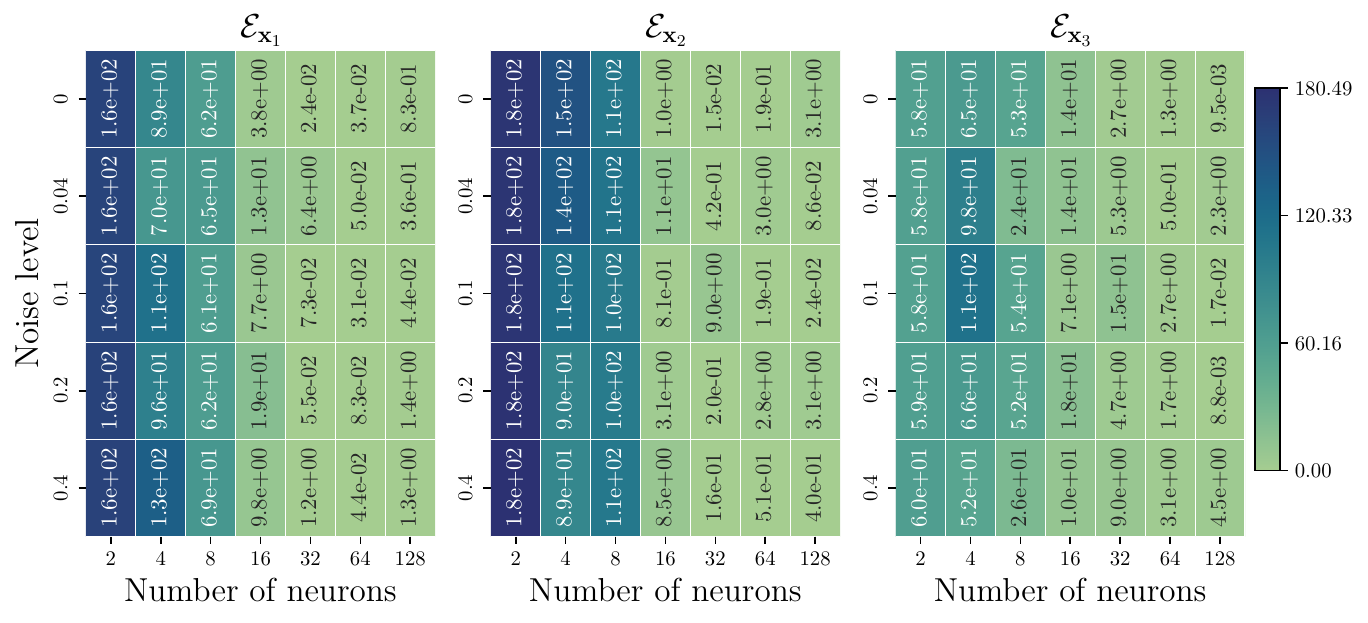}
		\caption{Using \deepymod.}
	\end{subfigure}
	\caption{Lorenz example:  A comparison of \ineuralsindy\ and \deepymod\ under \sceneb\ with the scaling factor $\alpha=0.1$.}
	\label{heatmap_lorenz_fixed_samples}	
\end{figure}

\begin{figure}[!tb]
	\centering
	
	\begin{subfigure}[b]{0.49\textwidth}
		\centering
		\includegraphics[width=1\textwidth]{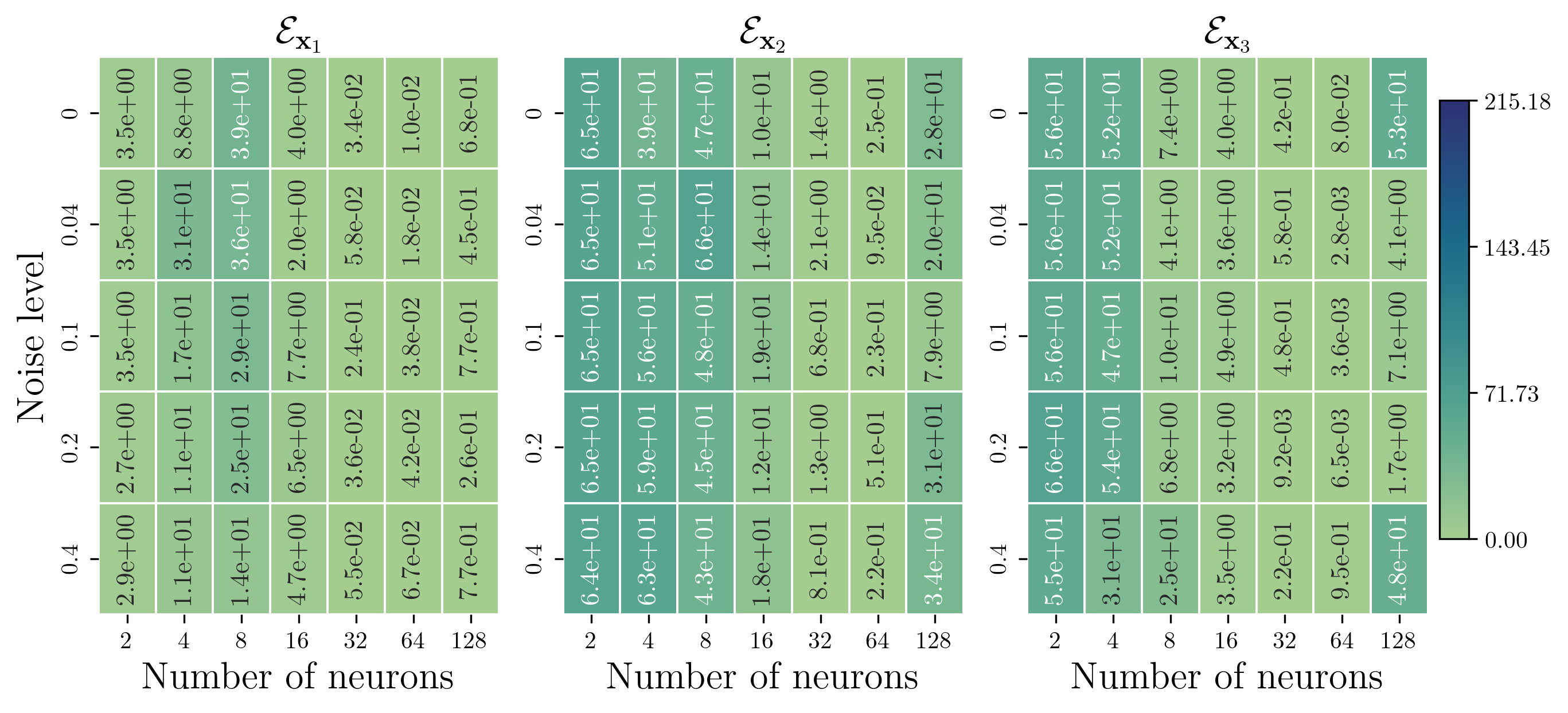}
		\caption{Using \ineuralsindy.}
	\end{subfigure}
	\begin{subfigure}[b]{0.49\textwidth}
		\centering
		\includegraphics[width=1\textwidth]{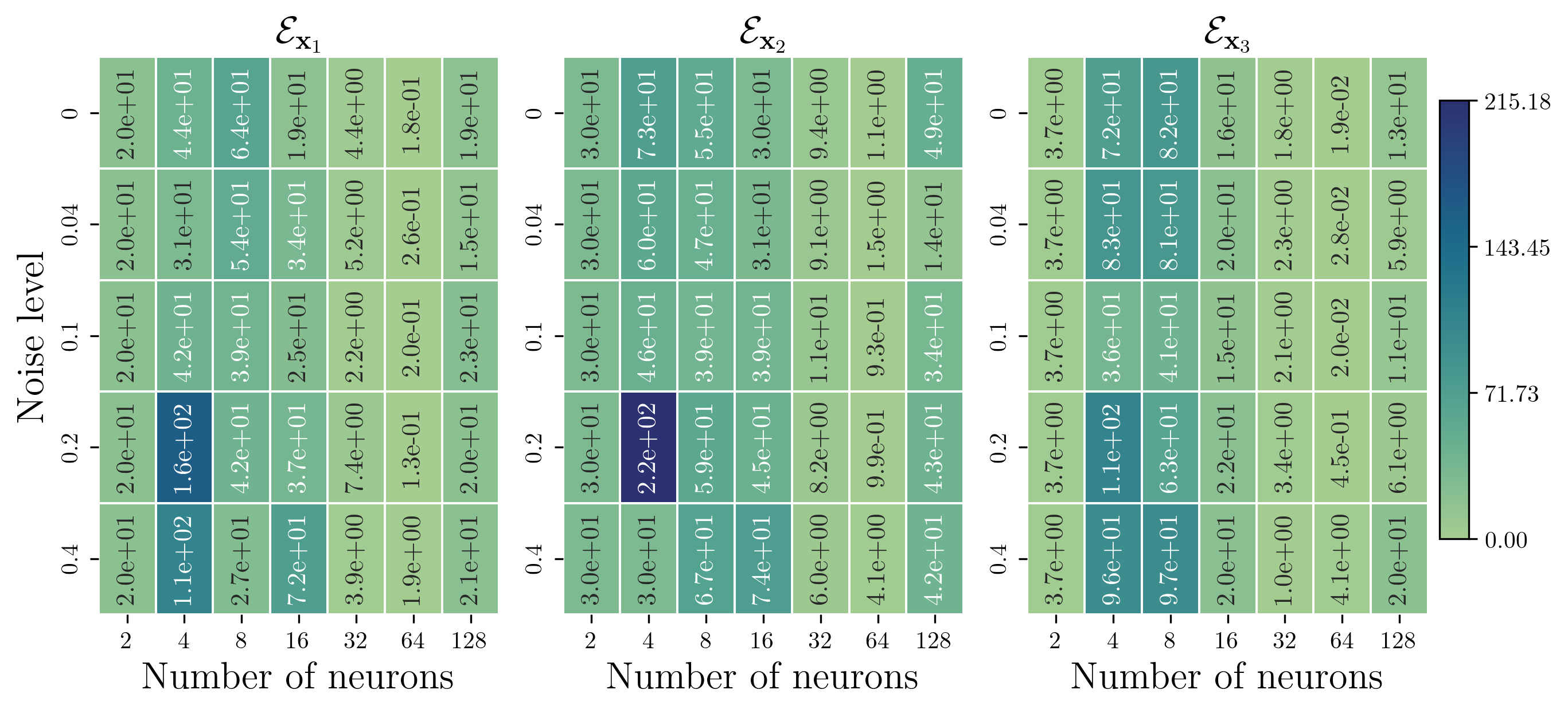}
		\caption{Using \deepymod.}
	\end{subfigure}
	\caption{Lorenz example: A comparison of \ineuralsindy\ and \deepymod\ under \sceneb\ with the scaling factor $\alpha=1$}
\label{heatmap_lorenz_fixed_samples_scaling_factor_1}		
\end{figure}

\paragraph{Results:} We conduct experiments using a scaling factor $\alpha = 0.1$.  Further,  we aim to learn governing equations from the noisy data with noise levels of $\sigma=\{0,\ 0.04,\ 0.1,\ 0.2,\ 0.4 \}$. We report the obtained results in \Cref{Lorenz_table}, where we notice that  \ineuralsindy\ and \deepymod\ yield similar performance except for the case of higher noise level (e.g., see the results for $\sigma = 0.4$), where \ineuralsindy\ recovered the equations better. However, \rksindy\ performs poorly for the higher noise levels. 
Next, we conduct a performance analysis of \ineuralsindy\ and \deepymod\ for \scenea\ and \sceneb. We note that in both scenarios,  the training data are generated using a single initial condition  $({\bx_1(0)}, {\bx_2(0)}, {\bx_3(0)})=(-8,\ 7,\ 27)$. 
\begin{itemize}
	\item \scenea: We compare \ineuralsindy\ and \deepymod\ under \scenea. We also investigate the effect of the scaling factor $\alpha$ and consider two values of it, i.e., $\alpha=\{0.1,\ 1\}$. We fix the \DNN\ architecture to have three hidden layers, each having $64$ neurons. We consider different sample sizes and noise levels to compare the performance of \ineuralsindy\ and \deepymod. For $\alpha = 0.1$, we show the results in  \Cref{heatmap_lorenz_fixed_neurons_sacling_0.1}, where we notice that \ineuralsindy\ outperforms \deepymod\ in most cases. A similar observation is made for $\alpha =1$, which is reported in 
	\Cref{heatmap_lorenz_fixed_neurons_sacling_1}. For these experiments, it is hard to conclude the effect of the scaling factors, as we notice that in some cases, the scaling improves the performance, and in some cases, it is not the case.
	
	\item \sceneb: In this case, we fix the sample size to $400$. We also fix the number of hidden layers for the \DNN\ architecture to three but vary the number of neurons in each layer. We also conduct experiments to see the effect of the scaling factor in this case as well. The results for $\alpha = 1$ and $\alpha = 0.1$ are shown in \Cref{heatmap_lorenz_fixed_samples} and \Cref{heatmap_lorenz_fixed_samples_scaling_factor_1}, respectively. These heat maps indicate the outperformance of \ineuralsindy~in most cases. We also observe that for a larger number of neurons, the scaling factor $\alpha=0.1$ slightly performs better compared to scaling factor $\alpha=1$ in both \ineuralsindy\ as well as \deepymod.  
\end{itemize}

\paragraph{A comparison of  \ineuralsindy\ with WEAK-SINDy:} Beside our previous comprehensive study, we next compare \ineuralsindy\ with \weaksindy, which also does not require any estimate of derivatives using noisy data; for more details on \weaksindy, we refer to \cite{Weak_SINDy}. 

For this study, we again consider the same initial condition as used for \scenea\ and \sceneb. We take  $2000$ data points in  the time interval $[0,10]$, which are corrupted using different noise levels  $\sigma = \{0, 0.02, 0.08, 0.1\}$. To discover the governing equations, we consider a dictionary of polynomials up to degree two. For training \ineuralsindy\, we use the same setting as discussed in \Cref{subsec:lorenzsetip}. For \weaksindy, we consider the code provided by the authors\footnote{\url{https://github.com/MathBioCU/WSINDy_ODE/tree/master}}. We report the results in \Cref{weaksindy_table}, which indicates that  \ineuralsindy\ outperforms \weaksindy~in the presence of high noise.

\section{Conclusions}\label{Conclusion}
In this work, we proposed a methodology, namely \ineuralsindy, to discover governing equations using noisy and scarce data. It consists of three main components---these are: (a) learning an implicit representation based on given noisy data using a deep neural network, (b) setting up a sparse regression problem inspired by \sindy\ \cite{brunton2016discovering}, discovering governing equations, and (c) utilize an integral form of differential equations. We have combined all these components innovatively to learn governing equations from noisy data. Particularly, we highlight that we leverage the implicit representation using neural networks to estimate the derivative using automatic differential to avoid any numerical derivative estimation using noisy data. 
We have shown how \ineuralsindy\ can be employed when data are collected using multiple trajectories. Furthermore, we have presented an extensive comparison of the proposed methodology with \rksindy\ \cite{Goyal_2022} and \deepymod\ \cite{Deepymod_2021}, where we noticed that \ineuralsindy\ clearly out-performed \rksindy, and in many cases, \ineuralsindy\ also yielded better or comparable results as compared to \deepymod, expect for the FHN example. We also compared \ineuralsindy\ with \weaksindy~using the Lorenz example, where we noticed a better performance of \ineuralsindy. 
In the future, we would like to extend the proposed framework to the identification of parametric and control-driven dynamical systems. We also like to combine the idea of the ensemble discussed in \cite{Fasel_2022_ensemble} to further improve the quality of learned governing equations. 

\addcontentsline{toc}{section}{References}
\bibliographystyle{IEEEtran}
\bibliography{mybib,iNeuralSINDyref}

\appendix
\section{Appendix}

In this appendix, we present a comparison of learned governing equations using various methods. Moreover, the code and examples are available at the following Git repository \footnote{\url{https://github.com/Ali-Forootani/iNeural_SINDy_paper/tree/main}}.

\begin{table}[h]
	\caption{Estimated coefficients for linear damped oscillator}
	\label{2_D_Oscillators_table}

	\centering
	\tiny
	\begin{tabular}{|p{0.5cm} |p{3.5cm}|p{3.6cm}|p{3.5cm}|}	
		\hline
		
		\multicolumn{4}{|c|}{Learned equations} \\
		\hline
		
		\centering Noise level & \centering \ineuralsindy &  \centering \deepymod &  \multicolumn{1}{c|}{\rksindy~\cite{Goyal_2022}}  
		
		\\
		\hline
		\vspace{0.02cm}
		\centering

		$0.00$
		&
		\vspace{0.02cm}
		
		$
		\begin{aligned}
			\dot{\bx}_1(t)&=-0.100\bx_1(t) + 2.000\bx_2(t)\nonumber\\
			\dot{\bx}_2(t)&= -2.000\bx_1(t) - 0.100\bx_2(t)\nonumber
		\end{aligned}
		$
		
		&
		\vspace{0.02cm}
		
		$
		\begin{aligned}
			\dot{\bx}_1(t)&=-0.100\bx_1(t) + 2.000\bx_2(t)\nonumber\\
			\dot{\bx}_2(t)&= -2.000\bx_1(t) - 0.100\bx_2(t)\nonumber
		\end{aligned} 
		$
		
		&
		\vspace{0.02cm}
		
		$
		\begin{aligned}
			\dot{\bx}_1(t)&=-0.100\bx_1(t) + 2.000\bx_2(t)\nonumber\\
			\dot{\bx}_2(t)&= -2.000\bx_1(t) - 0.100\bx_2(t)\nonumber
		\end{aligned}
		$
		
		\\
		\hline
		\vspace{0.02cm}
		\centering
		$0.02$
		&
		\vspace{0.02cm}
		$
		\begin{aligned}
			\dot{\bx}_1(t)&=-0.100 \bx_1(t)+ 2.000 \bx_2(t)\nonumber\\
			\dot{\bx}_2(t)&= -2.000\bx_1(t)-0.100 \bx_2(t)\nonumber
		\end{aligned}
		$
		\vspace{0.02cm}
		&
		\vspace{0.02cm}
		$
		\begin{aligned}
			\dot{\bx}_1(t)&=-0.100x_1(t) + 2.000\bx_2(t)\nonumber\\
			\dot{\bx}_2(t)&= -2.000\bx_1(t) - 0.099\bx_2(t)\nonumber
		\end{aligned}
		$
		\vspace{0.02cm}
		&
		\vspace{0.02cm}
		$
		\begin{aligned}
			\dot{\bx}_1(t)&=-0.103 \bx_1(t)+ 1.997 \bx_2(t)\nonumber\\
			\dot{\bx}_2(t)&= -2.003\bx_1(t)-0.104 \bx_2(t)\nonumber
		\end{aligned}
		$
		\vspace{0.02cm}
		\\
		\hline
		\vspace{0.02cm}
		\centering
		$0.04$
		&
		\vspace{0.02cm}
		$
		\begin{aligned}
			\dot{\bx}_1(t)&=-0.105 \bx_1(t)+  2.007\bx_2(t)\nonumber\\
			\dot{\bx}_2(t)&= -1.994\bx_1(t)- 0.094\bx_2(t)\nonumber
		\end{aligned}
		$
		\vspace{0.02cm}
		&
		\vspace{0.02cm}
		$
		\begin{aligned}
			\dot{\bx}_1(t)&=-0.098\bx_1(t) + 2.000\bx_2(t)\nonumber\\
			\dot{\bx}_2(t)&= -2.000\bx_1(t) - 0.099\bx_2(t)\nonumber
		\end{aligned}
		$
		\vspace{0.02cm}
		&
		\vspace{0.02cm}
		$
		\begin{aligned}
			\dot{\bx}_1(t)&= -0.109 \bx_1(t) + 1.989 \bx_2(t)\nonumber\\
			\dot{\bx}_2(t)&= -2.009 \bx_1(t) - 0.110 \bx_2(t)\nonumber
		\end{aligned}
		$
		\vspace{0.02cm}
		\\
		\hline
		\vspace{0.02cm}
		\centering
		$0.08$
		&
		\vspace{0.02cm}
		$
		\begin{aligned}
			\dot{\bx}_1(t)=-0.097 \bx_1(t)+  1.997 \bx_2(t)\nonumber\\
			\dot{\bx}_2(t)= -2.003\bx_1(t)- 0.106\bx_2(t)\nonumber
		\end{aligned}
		$
		\vspace{0.02cm}
		&
		\vspace{0.02cm}
		$
		\begin{aligned}
			\dot{\bx}_1(t)=-0.1000 \bx_1(t) + 1.996 \bx_2(t)\nonumber\\
			\dot{\bx}_2(t)=-2.005 \bx_1(t) - 0.102 \bx_2(t)\nonumber
		\end{aligned}
		$
		\vspace{0.02cm}
		&
		\vspace{0.02cm}
		$
		\begin{aligned}
			\dot{\bx}_1(t)= -0.170 \bx_1(t) + 1.916 \bx_2(t)\nonumber\\
			\dot{\bx}_2(t)= -2.073 \bx_1(t) - 0.177 \bx_2(t)\nonumber
		\end{aligned}
		$
		\vspace{0.02cm}
		\\
		\hline
		
	\end{tabular}
	
\end{table}

\normalsize


	\begin{table}[h]
		\centering
		\caption{Estimated coefficients for Cubic damped oscillator}
		\label{Cubic oscillator}
		
		\tiny
		
		\begin{tabular}
			{|p{0.5cm} |p{3.6cm}|p{3.6cm}|p{5.5cm}|}
			
			\hline
			\multicolumn{4}{|c|}{Estimated System} \\
			\hline
			\centering Noise level& \centering \ineuralsindy&  \centering \deepymod &  \multicolumn{1}{c|}{\rksindy~\cite{Goyal_2022}}  \\
			\hline
			\vspace{0.02cm}
			\centering
			$0.00$
			&
			
			\vspace{0.02cm}
			$
			\begin{aligned}
			\dot{\bx}_1(t)&= -0.100 \bx_1^3(t) + 2.000 \bx_2^3(t)\nonumber\\
			\dot{\bx}_2(t)&= -2.000\bx_1^3(t) - 0.100 \bx_2^3(t)\nonumber
			\end{aligned}
			$
			
			&
			
			\vspace{0.02cm}
			$
			\begin{aligned}
			\dot{\bx}_1(t)&= -0.100 \bx_1^3(t) + 2.000 \bx_2^3(t)\nonumber\\
			\dot{\bx}_2(t)&= -2.000\bx_1^3(t)- 0.100 \bx_2^3(t)\nonumber
			\end{aligned}
			$
			
			&
			\vspace{0.02cm}
			$
			\begin{aligned}
			\dot{\bx}_1(t)&= -0.099 \bx_1^3(t) + 2.000 \bx_2^3(t)\nonumber\\
			\dot{\bx}_2(t)&= -2.000\bx_1^3(t)- 0.100 \bx_2^3(t)\nonumber
			\end{aligned}
			$
			
			\\
			\hline
			\vspace{0.02cm}
			\centering
			$0.02$
			&
			\vspace{0.02cm}
			$
			\begin{aligned}
			\dot{\bx}_1(t)&= -0.103 \bx_1^3(t) + 1.996 \bx_2^3(t)\nonumber\\
			\dot{\bx}_2(t)&= -1.997\bx_1^3(t)- 0.099 \bx_2^3(t)\nonumber
			\end{aligned}
			$
			\vspace{0.02cm}
			&
			\vspace{0.02cm}
			$
			\begin{aligned}
			\dot{\bx}_1(t)& = -0.098 \bx_1^3(t) + 2.005 \bx_2^3(t)\nonumber\\
			\dot{\bx}_2(t)& = -1.995\bx_1^3(t)- 0.102 \bx_2^3(t)\nonumber
			\end{aligned}
			$
			\vspace{0.02cm}
			&
			\vspace{0.02cm}
			$
			\begin{aligned}
			\dot{\bx}_1(t)&=-0.102 \bx_1^3(t)+ 2.004 \bx_2^3(t)\nonumber\\
			\dot{\bx}_2(t)&= -1.987\bx_1(t) -0.054\bx_1(t)\bx_2^2(t) -0.119 \bx_2^3(t)\nonumber
			\end{aligned}
			$
			\vspace{0.02cm}
			\\
			\hline
			\vspace{0.02cm}
			\centering
			$0.04$
			&
			\vspace{0.02cm}
			$
			\begin{aligned}
			\dot{\bx}_1(t)& = -0.105 \bx_1^3(t) + 1.994 \bx_2^3(t)\nonumber\\
			\dot{\bx}_2(t)& = -1.989\bx_1^3(t)- 0.0984 \bx_2^3(t)\nonumber
			\end{aligned}
			$
			\vspace{0.02cm}
			&
			\vspace{0.02cm}
			$
			\begin{aligned}
			\dot{\bx}_1(t)&= -0.0974 \bx_1^3(t) + 2.000 \bx_2^3(t)\nonumber\\
			\dot{\bx}_2(t)&= -2.015 \bx_1^3(t)- 0.104 \bx_2^3(t)\nonumber
			\end{aligned}
			$
			\vspace{0.02cm}
			&
			\vspace{0.02cm}
			$
			\begin{aligned}
			\dot{\bx}_1(t) &= 0.076\bx_2(t)-0.147 \bx_1^3(t) \\&+ 0.059 \bx_1(t) \bx_2^2(t)  + 1.907 \bx_2^3(t)\nonumber\\
			\dot{\bx}_2(t)&= -0.051\bx_1^2(t)- 0.059 \bx_1(t)\bx_2(t) \\&-2.042\bx_1^3(t) + 0.112 \bx_1^2(t) \bx_2(t) - 0.167 \bx_2^3(t)\nonumber
			\end{aligned}
			$
			\vspace{0.02cm}
			\\
			\hline
			\vspace{0.02cm}
			\centering
			$0.06$
			\vspace{0.02cm}
			&
			\vspace{0.02cm}
			$
			\begin{aligned}
			\dot{\bx}_1(t)&=-0.100 \bx_1^3(t)+  1.947\bx_2^3(t)\nonumber\\
			\dot{\bx}_2(t)&= -1.986\bx_1^3(t)- 0.112\bx_2^3(t)\nonumber
			\end{aligned}
			$
			\vspace{0.02cm}
			&
			\vspace{0.02cm}
			$
			\begin{aligned}
			\dot{\bx}_1(t) &=-0.097 \bx_1^3(t)+  1.972 \bx_2^3(t)\nonumber\\
			\dot{\bx}_2(t) &= -2.040 \bx_1^3(t)- 0.105 \bx_2^3(t)\nonumber
			\end{aligned}
			$
			\vspace{0.02cm}
			&
			\vspace{0.02cm}
			
			$
			\begin{aligned}
			\dot{\bx}(t)&= 0.156 \bx_1(t) + 0.167 \bx_2(t) - 0.071 \bx_2^2(t) \\&- 0.294 \bx_1^3(t) - 0.165 \bx_1(t)\bx_2^2(t) + 1.688 \bx_2^3(t) \\
			\dot{\bx}_2(t)&= 0.107 \bx_1(t)- 0.181 \bx_2(t) +  0.084 \bx_1(t) \bx_2(t) \\& - 0.084 \bx_2^2(t) - 2.169\bx_1^3(t) -0.134 \bx_1(t) \bx_2^2(t) \\&-0.355 \bx_2^3(t)\nonumber
			\end{aligned}
			$
			\vspace{0.02cm}
			\\
			\hline
		\end{tabular}
	\end{table}


\begin{table}[h]
	\centering
	\caption{Estimated coefficients for Fitz-Hugh Nagumo}
	\label{Fitz-Hugh_Nagumo_table}
	
	\tiny
	
	\begin{tabular}{|p{0.5cm} |p{3.5cm}|p{3.5cm}|p{4cm}|  }
		\hline
		\multicolumn{4}{|c|}{Estimated System} \\
		\hline
		\centering Noise level& \centering \ineuralsindy& \centering \deepymod & \multicolumn{1}{c|}{\rksindy~\cite{Goyal_2022}} \\
		\hline
		\vspace{0.02cm}
		\centering
		$
		0.00
		$
		&
		\vspace{0.02cm}
		$
		\begin{aligned}
		\dot{\bx}_1(t) &= 0.989 \bx_1(t) - 0.993 \bx_2(t) \\ &- 0.329 \bx_1^3(t) + 0.100\nonumber\\
		\dot{\bx}_2(t) &= 0.100 \bx_1(t) - 0.099\bx_2(t)\nonumber
		\end{aligned}
		$
		
		&
		\vspace{0.02cm}
		$
		\begin{aligned}
		\dot{\bx}_1(t) &= 0.992 \bx_1(t) - 0.994 \bx_2(t) \\ &- 0.330 \bx_1^3(t) + 0.100\nonumber\\
		\dot{\bx}_2(t) &= 0.100 \bx_1(t) - 0.100\bx_2(t)\nonumber
		\end{aligned}
		$
		
		&
		\vspace{0.02cm}
		$
		\begin{aligned}
		\dot{\bx}_1(t) &= 0.986 \bx_1(t) - 0.996 \bx_2(t) \\ &- 0.328 \bx^3_1(t) + 0.0993\nonumber\\
		\dot{\bx}_2(t) &= 0.100 \bx_1(t) - 0.099\bx_2(t)\nonumber
		\end{aligned}
		$
		
		\\	
		\hline
		\vspace{0.02cm}
		\centering
		$0.02$
		&
		\vspace{0.02cm}
		$
		\begin{aligned}
		\dot{\bx}_1(t) &= 0.993 \bx_1(t) - 0.997 \bx_2(t) \\ &- 0.330 \bx^3_1(t) + 0.100\nonumber\\
		\dot{\bx}_2(t) &= 0.100 \bx_1(t) - 0.100\bx_2(t)\nonumber
		\end{aligned}
		$
		
		&
		\vspace{0.02cm}
		$
		\begin{aligned}
		\dot{\bx}_1(t) &= 0.994 \bx_1(t) - 0.996 \bx_2(t) \\ &- 0.330 \bx_1(t)^3 + 0.100\nonumber\\
		\dot{\bx}_2(t) &= 0.100 \bx_1(t) - 0.100\bx_2(t)\nonumber
		\end{aligned}
		$
		
		&
		\vspace{0.02cm}
		$
		\begin{aligned}
		\dot{\bx}_1(t) &= 1.000 \bx_1(t) - 1.000 \bx_2(t) \\ &- 0.333 \bx_1(t)^3 + 0.100\nonumber\\
		\dot{\bx}_2(t) &= 0.100 \bx_1(t) - 0.101\bx_2(t)\nonumber
		\end{aligned}
		$
		
		\\
		\hline
		\vspace{0.02cm}
		\centering
		$0.04$
		&
		\vspace{0.02cm}
		$
		\begin{aligned}
		\dot{\bx}_1(t) &= 0.997 \bx_1(t) - 0.999 \bx_2(t) \\ &- 0.332 \bx_1(t)^3 + 0.101\nonumber\\
		\dot{\bx}_2(t) &= 0.100 \bx_1(t) - 0.102\bx_2(t)\nonumber
		\end{aligned}
		$
		\vspace{0.02cm}
		&
		\vspace{0.02cm}
		$
		\begin{aligned}
		\dot{\bx}_1(t) &= 0.967 \bx_1(t) - 0.973 \bx_2(t) \\ &- 0.321 \bx_1(t)^3 + 0.098\nonumber\\
		\dot{\bx}_2(t) &= 0.100 \bx_1(t) - 0.100 \bx_2(t)\nonumber
		\end{aligned}
		$
		
		&
		\vspace{0.02cm}
		$
		\begin{aligned}
		\dot{\bx}_1(t) &= 1.000 \bx_1(t) - 1.010 \bx_2(t) \\ &- 0.335  \bx_1(t)^3 + 0.1\nonumber\\
		\dot{\bx}_2(t) &= 0.101 \bx_1(t) - 0.109 \bx_2(t)\nonumber
		\end{aligned}
		$
		
		\\
		\hline
		\vspace{0.02cm}
		\centering
		$0.08$
		&
		\vspace{0.02cm}
		$
		\begin{aligned}
		\dot{\bx}_1(t) &= 0.978 \bx_1(t) - 0.985 \bx_2(t) \\ &- 0.324 \bx_1(t)^3 + 0.097  \nonumber\\
		\dot{\bx}_2(t) &= 0.100 \bx_1(t) - 0.112 \bx_2(t)\nonumber
		\end{aligned}
		$
		
		&
		\vspace{0.02cm}
		$
		\begin{aligned}
		\dot{\bx}_1(t) &= 0.908 \bx_1(t) - 0.941 \bx_2(t) \\ &- 0.297 \bx_1(t)^3 +  0.096\nonumber\\
		\dot{\bx}_2(t) &= 0.102 \bx_1(t) - 0.106 \bx_2(t)\nonumber
		\end{aligned}
		$
		
		&
		\vspace{0.02cm}
		$
		\begin{aligned}
		\dot{\bx}_1(t) & = 0.336 \bx_1(t) - 0.081 \bx^2_2(t) \\ & -0.076 \bx^3_1 -0.117 \bx_1^2(t) \bx_2(t)\\ &-0.139 \bx_1(t) \bx^2_2(t) - 0.113 \bx^3_2(t) \nonumber\\
		\dot{\bx}_2(t) &= 0.098 \bx_1(t) - 0.062 \bx_2(t)^2 \\ &- 0.220 \bx_2^3\nonumber
		\end{aligned}
		$
		
		\\
		\hline	
	\end{tabular}
\end{table}


\begin{table}[h]
	\centering
	\caption{Estimated coefficients for Lorenz}
	\label{Lorenz_table}
	\tiny
	
	\begin{tabular}{|p{0.5cm}|p{4.1cm}|p{4.1cm}|p{4.1cm}|}
		\hline
		\multicolumn{4}{|c|}{Estimated System} \\
		\hline
		\centering Noise level& \centering \ineuralsindy& \centering \deepymod & \multicolumn{1}{c|}{\rksindy~\cite{Goyal_2022}} \\
		\hline
		\vspace{0.02cm}
		\centering
		$0.00$
		&
		\vspace{0.02cm}
		$
		\begin{aligned}
		\dot{\bx}_1(t) & =  -9.989 \bx_1(t)+ 9.991 \bx_2(t)\nonumber\\
		\dot{\bx}_2(t) & = 28.022 \bx_1(t) - 1.004 \bx_2(t) \\&+ 10.006 \bx_1(t) x_3(t) \nonumber\\
		\dot{\bx}_3(t) & = -2.666 \bx_3(t) + 10.013 \bx_1(t) x_2(t)
		\end{aligned}
		$
		
		&
		\vspace{0.02cm}
		$
		\begin{aligned}
		\dot{\bx}_1(t) & =  -9.992 \bx_1(t) + 9.993 \bx_2(t) \nonumber\\
		\dot{\bx}_2(t) & = 28.026 \bx_1(t) - 1.000 \bx_2(t) \\& + 10.000 \bx_1(t) \bx_3(t) \nonumber\\
		\dot{\bx}_3(t) & = -2.666 \bx_3(t) + 10.000 \bx_1(t) \bx_2(t)
		\end{aligned}
		$
		
		&
		\vspace{0.02cm}
		$
		\begin{aligned}
		\dot{\bx}_1(t) & =  -9.995 \bx_1(t) + 10.003 \bx_2(t) \nonumber\\
		\dot{\bx}_2(t) & = 28.032 \bx_1(t) - 1.000 \bx_2(t) \\&+ 10.014 \bx_1(t) \bx_3(t) \nonumber\\
		\dot{\bx}_3(t) & = -2.665 \bx_3(t) + 10.020 \bx_1(t) \bx_2(t)
		\end{aligned}
		$

		\\
		\hline
		\vspace{0.02cm}
		\centering			
		$0.10$
		&
		\vspace{0.02cm}
		$
		\begin{aligned}
		\dot{\bx}_1(t) & = -9.982 \bx_1(t) + 9.985 \bx_2(t)\nonumber\\
		\dot{\bx}_2(t) & = 28.001 \bx_1(t) - 1.001 \bx_2(t) \\&- 10.005 \bx_1(t)\bx_3(t) \nonumber\\
		\dot{\bx}_3(t) & = -2.666 \bx_3(t) + 10.008 \bx_1(t) x_2(t)  
		\end{aligned}
		$
		
		&
		\vspace{0.02cm}
		$
		\begin{aligned}
		\dot{\bx}_1(t) & = -9.997 \bx_1(t) + 9.994 \bx_2(t)\nonumber\\
		\dot{\bx}_2(t) & = 28.028 \bx_1(t) - 1.005 \bx_2(t) \\&- 10.008 \bx_1(t)\bx_3(t) \nonumber\\
		\dot{\bx}_3(t) & = -2.666 \bx_3(t) + 9.999 \bx_1(t) \bx_2(t)  
		\end{aligned}
		$
		
		&
		\vspace{0.02cm}
		$
		\begin{aligned}
		\dot{\bx}_1(t) & = -10.004 \bx_1(t) + 10.008 \bx_2(t)\nonumber\\
		\dot{\bx}_2(t) & = 28.060 \bx_1(t) - 0.997 \bx_2(t) \\&- 10.025 \bx_1(t)\bx_3(t) \nonumber\\
		\dot{\bx}_3(t) & = -2.664 \bx_3(t) + 10.030 \bx_1(t) \bx_2(t)  
		\end{aligned}
		$
		
		\\
		\hline
		\vspace{0.02cm}
		\centering
		$0.20$
		&
		\vspace{0.02cm}
		$
		\begin{aligned}
		\dot{\bx}_1(t) & = -10.016 \bx_1(t) + 10.008 \bx_2(t)\nonumber\\
		\dot{\bx}_2(t) & = 27.980 \bx_1(t) - 0.991 \bx_2(t) \\&- 9.998 \bx_1(t)\bx_3(t) \nonumber\\
		\dot{\bx}_3(t) & = -2.670 \bx_3(t) + 10.010 \bx_1(t) \bx_2(t)  
		\end{aligned}
		$
		
		&
		\vspace{0.02cm}
		$
		\begin{aligned}
		\dot{\bx}_1(t) & = -10.034 \bx_1(t) + 10.032 \bx_2(t)\nonumber\\
		\dot{\bx}_2(t) & = 27.968 \bx_1(t) - 0.998 \bx_2(t) \\&- 9.986 \bx_1(t)\bx_3(t) \nonumber\\
		\dot{\bx}_3(t) & = -2.669 \bx_3(t) + 9.993 \bx_1(t) \bx_2(t)  
		\end{aligned}
		$
		
		&
		\vspace{0.02cm}
		$
		\begin{aligned}
		\dot{\bx}_1(t) & = -10.076 \bx_1(t) + 10.068 \bx_2(t)\nonumber\\
		\dot{\bx}_2(t) & = 27.846 \bx_1(t) - 0.928 \bx_2(t) \\&- 9.967 \bx_1(t)\bx_3(t) \nonumber\\
		\dot{\bx}_3(t) & = -2.674 \bx_3(t) + 9.995 \bx_1(t) \bx_2(t)  
		\end{aligned}
		$
		
		\\
		\hline
		\vspace{0.02cm}
		\centering
		$0.40$
		&
		\vspace{0.02cm}
		$
		\begin{aligned}
		\dot{\bx}_1(t) & = -10.122 \bx_1(t) + 10.080 \bx_2(t)\nonumber\\
		\dot{\bx}_2(t) & = 27.714 \bx_1(t) - 0.920 \bx_2(t) \\&- 9.938 \bx_1(t)\bx_3(t) \nonumber\\
		\dot{\bx}_3(t) & = -2.673 \bx_3(t) + 9.933 \bx_1(t) \bx_2(t)  
		\end{aligned}
		$
		
		&	
		\vspace{0.02cm}
		$
		\begin{aligned}
		\dot{\bx}_1(t) & = -9.887 \bx_1(t) + 9.902 \bx_2(t)\nonumber\\
		\dot{\bx}_2(t) & = 26.627 \bx_1(t) - 9.596 \bx_1(t)\bx_3(t) \\&-0.3323 \bx_2(t) \bx_3(t) \nonumber\\
		\dot{\bx}_3(t) & = -2.657 \bx_3(t) + 0.035 \bx_1(t) \bx_3(t) \\&+ 9.964 \bx_1(t) \bx_2(t)  
		\end{aligned}
		\vspace{0.02cm}
		$
		
		&			
		\vspace{0.02cm}
		$
		\begin{aligned}
		\dot{\bx}_1(t) & = -10.976 \bx_1(t) + 10.272 \bx_2(t) \\&+ 0.225 \bx_1(t) \bx_3(t)\nonumber\\
		\dot{\bx}_2(t) & = 26.144 \bx_1(t) - 9.490 \bx_1(t)\bx_3(t) \\&-0.270 \bx_2(t) \bx_3(t) \nonumber\\
		\dot{\bx}_3(t) & = -2.683 \bx_3(t) + 10.030 \bx_1(t) \bx_2(t)  
		\end{aligned}
		$
		
		\\
		\hline		
	\end{tabular}
\end{table}



\begin{table}[h]
	\centering
	\caption{Comparison of \ineuralsindy\ with WEAK-SINDy}
	\label{weaksindy_table}
	
	\tiny
	\begin{tabular}{|p{0.5cm}|p{5.5cm}|p{6cm}|}
		\hline
		\multicolumn{3}{|c|}{Estimated System} \\
		\hline
		\centering Noise level& \centering \ineuralsindy & \multicolumn{1}{c|}{\weaksindy}\\
		\hline
		\vspace{0.02cm}
		\centering
		$0.00$
		&
		\vspace{0.02cm}
		$
		\begin{aligned}
		\dot{\bx}_1(t) & =  -9.990 \bx_1(t) + 9.990 \bx_2(t)\nonumber\\
		\dot{\bx}_2(t) & = 28.021 \bx_1(t) - 1.003 \bx_2(t) - 1.006 \bx_1(t) x_3(t) \nonumber\\
		\dot{\bx}_3(t) & = -2.666 \bx_3(t) + 1.001 \bx_1(t) x_2(t)
		\end{aligned}
		$

		&
		\vspace{0.02cm}
		$
		\begin{aligned}
		\dot{\bx}_1(t) & =  10.000 \bx_1(t) - 10.000 \bx_2(t) \nonumber\\
		\dot{\bx}_2(t) & = 28.026 \bx_1(t) - 1.000 \bx_2(t)  - 1.000 \bx_1(t) \bx_3(t) \nonumber\\
		\dot{\bx}_3(t) & = -2.670 \bx_3(t) + 1.000 \bx_1(t) \bx_2(t)
		\end{aligned}
		$
		\\
		\hline
		\vspace{0.02cm}
		\centering			
		$0.02$
		&
		\vspace{0.02cm}
		$
		\begin{aligned}
		\dot{\bx}_1(t) & = -9.989 \bx_1(t) + 9.991 \bx_2(t)\nonumber\\
		\dot{\bx}_2(t) & = 28.013 \bx_1(t) - 1.001 \bx_2(t) - 1.000 \bx_1(t)\bx_3(t) \nonumber\\
		\dot{\bx}_3(t) & = -2.666 \bx_3(t) + 1.001 \bx_1(t) x_2(t)  
		\end{aligned}
		$
		
		&
		\vspace{0.02cm}
		$
		\begin{aligned}
		\dot{\bx}_1(t) & = -9.908 \bx_1(t) + 9.926 \bx_2(t)\nonumber\\
		\dot{\bx}_2(t) & = 27.857 \bx_1(t) - 0.991 \bx_2(t) - 0.998 \bx_1(t)\bx_3(t) \nonumber\\
		\dot{\bx}_3(t) & = -2.658 \bx_3(t) + 1.002 \bx_1(t) \bx_2(t)  
		\end{aligned}
		$
		
		\\
		\hline
		\vspace{0.02cm}
		\centering
		$0.08$
		&

		\vspace{0.02cm}
		$
		\begin{aligned}
		\dot{\bx}_1(t) & = -9.989 \bx_1(t) + 9.991 \bx_2(t)\nonumber\\
		\dot{\bx}_2(t) & = 27.991 \bx_1(t) - 0.996 \bx_2(t) - 0.999 \bx_1(t)\bx_3(t) \nonumber\\
		\dot{\bx}_3(t) & = -2.666 \bx_3(t) + 1.001 \bx_1(t) \bx_2(t)  
		\end{aligned}
		$
		
		&
		\vspace{0.02cm}
		$
		\begin{aligned}
		\dot{\bx}_1(t) & = -9.798 \bx_1(t) + 9.778 \bx_2(t)\nonumber\\
		\dot{\bx}_2(t) & = 26.699 \bx_1(t) - 0.790 \bx_2(t) - 0.969 \bx_1(t)\bx_3(t) \nonumber\\
		\dot{\bx}_3(t) & = -2.640 \bx_3(t) + 1.010 \bx_1(t) \bx_2(t)  
		\end{aligned}
		$
		
		\\
		\hline
		\vspace{0.02cm}
		\centering
		$0.10$

		&
		\vspace{0.02cm}
		$
		\begin{aligned}
		\dot{\bx}_1(t) & = -9.988 \bx_1(t) + 9.991 \bx_2(t)\nonumber\\
		\dot{\bx}_2(t) & = 27.982 \bx_1(t) - 0.993 \bx_2(t) - 0.994 \bx_1(t)\bx_3(t) \nonumber\\
		\dot{\bx}_3(t) & = -2.666 \bx_3(t) + 1.019 \bx_1(t) \bx_2(t)  
		\end{aligned}
		$

		&	
		\vspace{0.02cm}
		$
		\begin{aligned}
		\dot{\bx}_1(t) & = -9.751 \bx_1(t) + 9.734 \bx_2(t)\nonumber\\
		\dot{\bx}_2(t) & = 26.337 \bx_1(t) - 0.716 \bx_1(t)\bx_3(t) -0.961 \bx_2(t) \bx_3(t) \nonumber\\
		\dot{\bx}_3(t) & = -2.630 \bx_3(t) + 1.010 \bx_1(t) \bx_2(t)  
		\end{aligned}
		\vspace{0.02cm}
		$
		
		\\
		\hline		
	\end{tabular}
\end{table}

%

\end{document}